\documentclass[journal]{IEEEtran}   
\usepackage{enumerate}
\usepackage{times}
\usepackage{url}
\usepackage[final]{graphicx}
\usepackage[reqno]{amsmath}
\usepackage{amsfonts}
\usepackage{multicol}
\usepackage{tabu}
\usepackage{times,amsmath,epsfig}
\usepackage{latexsym,amssymb}
\usepackage{cite}
\usepackage{subfigure}
\usepackage{amsmath}
\usepackage{amssymb}
\usepackage{array}
\usepackage[noend]{algorithmic}
\usepackage{algorithm}
\usepackage{color}
\usepackage{psfrag}
\usepackage{epsfig}

\newtheorem{pro}{Proposition}
\newtheorem{lemma}{Lemma}

\newtheorem{theorem}{Theorem}

\newcommand{\mb}{\mathbf}

\newcommand{\ba}{\mathbf{a}}

\newcommand{\bff}{\mathbf{f}}

\newcommand{\bl}{\mathbf{l}}

\newcommand{\bq}{\mathbf{q}}
\newcommand{\br}{\mathbf{r}}
\newcommand{\bs}{\mathbf{s}}

\newcommand{\bu}{\mathbf{u}}

\newcommand{\E}{\mathbb{E}}

\newcommand{\st}{\mbox{s.t.}}

\newcommand{\define}{{\triangleq}}
\newcommand{\argmin}{\arg\!\min}

\newcommand{\opt}{\textrm{opt}}
\newcommand{\ie}{i.e.}

\ifx\BlackBox\undefined
\newcommand{\BlackBox}{\rule{1.5ex}{1.5ex}}  
\fi

\ifx\proof\undefined
\newenvironment{proof}{\par\noindent{\em Proof:\ }}{\hfill\BlackBox\\[.0mm]}
\fi

\graphicspath{{./Figures/}}

\allowdisplaybreaks

\begin{document}

\title{
Phase Balancing Using Energy Storage\\ in Power Grids under Uncertainty
}

\author{Sun Sun, \IEEEmembership{Student Member,~IEEE},   Ben Liang, \IEEEmembership{Senior Member,~IEEE}, Min Dong, \IEEEmembership{Senior Member,~IEEE},   and Joshua A. Taylor, \IEEEmembership{Member,~IEEE}
\thanks{Sun Sun, Ben Liang, and Joshua A. Taylor are with the Department of  Electrical and Computer
Engineering, University of Toronto, Toronto, Ontario  M5S 3G4 (email: \{ssun,
liang\}@comm.utoronto.ca, josh.taylor@utoronto.ca).}
\thanks{Min Dong is with the Department of Electrical, Computer and Software Engineering, University of Ontario Institute of Technology, Oshawa, Ontario  L1H 7K4  (email:  min.dong@uoit.ca).}
\thanks{This work was supported by  the Natural Sciences and Engineering Research Council (NSERC) of Canada under Discovery Grants RGPIN-2014-05181 and RGPIN-2015-05506.}
\vspace{-1cm}}

\maketitle

\begin{abstract}
Phase balancing is essential to safe  power system operation. 
We consider a substation connected to multiple phases, each with single-phase loads, generation, and energy storage. A representative of the substation operates the system and aims to minimize the  cost of all phases and to balance loads among phases. We first consider ideal energy storage with lossless charging and discharging, and propose both centralized and distributed real-time algorithms taking into account  system uncertainty. The proposed algorithm does not require any system statistics and asymptotically achieves the minimum system cost with large energy storage. We then extend the algorithm to accommodate more realistic  non-ideal energy storage that has imperfect charging and discharging. The  performance of the proposed algorithm is evaluated through extensive simulation and compared with that of a benchmark greedy algorithm. Simulation shows that our algorithm leads to  strong performance over a wide range of storage characteristics. 

\end{abstract}


\section{Introduction}
In North America, many residential customers are connected to distribution systems through single-phase  lines. Phase balancing, \ie, maintaining the balance of  loads among phases, is crucial for power grid operation\cite{bkeps}. This is because phase imbalance  can increase energy losses and the risk of failures, and can also  degrade system power quality. With the spread of single-phase renewable generators, such as wind and solar generators, and large loads, such as electric vehicles,  phase imbalance  could be aggravated and thus deserves more careful study.  For example,  
the impact of  integration of electric vehicles on phase imbalance was investigated  in \cite{psj09}.

Previous works on phase balancing have considered methods such as phase swapping (e.g.,\cite{zcz98}) and feeder reconfiguration (e.g., \cite{hhlc93}). However, these  approaches can be ineffective or can incur extra costs on human resources,  maintenance expenses, and planned outage duration\cite{zcz98}. An alternative method is to employ energy storage to mitigate the imbalance among phases, which is the focus of this paper. 

Energy storage has been  used widely in power grids for  applications such as energy arbitrage, regulation, and load following \cite{dek10}. Examples of single-phase storage include:
\begin{itemize}
	\item Traditional standalone storage such as batteries, flywheels, etc\cite{cg14}.
	\item Batteries in single-phase connected buildings such as plug-in electric vehicles \cite{gwa11}.
	\item Aggregations of small single-phase deferrable loads, e.g., residential thermostatically controlled loads or electric vehicle garages, which have been shown to be representable as equivalent storage \cite{hsp13,hspv14,nts13}. 
\end{itemize}

The control of energy storage for power grid applications is, however, generally a challenging problem due to storage characteristics as well as  system uncertainty.
There are many existing works on storage control in power grids.
For example, using stochastic dynamic programming, the authors of \cite{sg13} proposed a stationary optimal policy for power balancing, and the authors of \cite{tcp13} investigated both optimal and suboptimal polices for energy balancing. Nevertheless, the derivation of an optimal policy under dynamic programming generally  relies on  system statistics and some specific form of the problem  structure, and therefore it cannot be easily extended. Similarly,  the authors of \cite{zh14} considered stochastic model predictive control. However, the algorithm performance can only be evaluated through numerical examples. 

Besides the above two approaches,  several recent works have employed  Lyapunov optimization \cite{bkneely} for energy storage control. 
In particular, the authors of \cite{rah11} investigated a power-cost minimization problem in data centers with energy storage, and were the first to use the technique of Lyapunov optimization for real-time storage control. The technique was then employed in several subsequent works to design energy storage control for various applications in  grid operation, such as power balancing \cite{sun14jsp,sun14tsg}, demand side management \cite{hwr12,gpfk13}, and EV charging \cite{lcnf11}. Furthermore, the authors of \cite{lqp14} analyzed the trade-off between averaging out the energy fluctuation across time and across space, the authors of \cite{qcyr14} studied generalized storage control with  general cost functions, and the authors of \cite{qcyr14b} investigated the management of networked storage with a DC power flow model.  Among these works, single storage control was considered in \cite{rah11,hwr12}, and \cite{qcyr14}. For multiple storage control,  charging efficiency  was incorporated into the storage model in \cite{lcnf11},  both charging and discharging efficiencies were introduced in \cite{sun14jsp},  and storage efficiency that models the energy loss over time was included in \cite{gpfk13}.
 
In this paper, we study the problem of phase balancing with energy storage  in the presence of  system uncertainty.  
Unlike prior works such as \cite{ams15} and \cite{bsl15}  that focus on heuristic algorithms for storage control in phase balancing, in this paper, we provide efficient algorithms with strong theoretical performance guarantee.
We consider a substation connected to multiple phases, each with single-phase uncontrollable flow, controllable flow, and  energy storage. In particular, we consider phase balancing on a time scale of  seconds to minutes. As such, we do not model power system physics such as frequency and voltage magnitude. 
Aiming at minimizing  the cost of all phases and mitigating phase imbalance, 
we propose  a real-time algorithm that can be easily implemented by the substation. Moreover, for the likely scenario  of limited communication between the substation and each phase, we provide a  distributed implementation of the real-time algorithm where only limited information exchange is required.

The main contributions of this paper are summarized as follows. First, we formulate a stochastic optimization problem for phase balancing incorporating   system uncertainty, storage characteristics, and power network constraints. Second, for ideal energy storage with lossless charging and discharging, we provide a real-time algorithm building on the framework of  Lyapunov optimization  and prove its analytical performance guarantee. Moreover, we offer distributed implementation of the algorithm with fast convergence. Third, we extend the algorithm to accommodate  non-ideal energy storage with imperfect charging and discharging efficiency and show its analytical performance. Finally,  to numerically evaluate the performance of the proposed algorithm, we compare it with a benchmark greedy algorithm under various settings and parameters. Simulation reveals that our proposed algorithm is competitive in general. In particular, the proposed algorithm has strong performance when applied to storage with a large energy capacity, a high value of the energy-power ratio (e.g., compressed air energy storage and batteries),  and moderate-to-high charging and discharging efficiency (e.g., the round-trip efficiency of storage is greater than $65\%$). In addition, a practical outcome of our analysis shows the following design guideline:  optimal power balancing favors even allocation of storage capacity over the phases.

Our paper is  technically most similar to \cite{qcyr14b}, in which a distributed real-time algorithm is  proposed for power grids with energy storage. However, these two papers are different in terms of the application, objective, communication topology, and power network constraints. Hence, the problem formulation and the design of distributed implementation are largely different. Moreover, for analysis,   charging and discharging efficiencies were not considered in the storage model in \cite{qcyr14b}. While the authors stated that their framework could further incorporate imperfect charging and discharging, no implementable algorithm was given to address that. In contrast, in this paper, we  provide an efficient algorithm to deal with imperfect charging and discharging.

A preliminary version of this work has been presented in \cite{sun15acc}. In this paper, we significantly extend \cite{sun15acc} in two ways: analytically, for the proposed algorithms, we provide more in-depth performance analysis for both ideal and non-ideal storage; numerically,  we implement extensive simulation by examining various storage characteristics and the effect of correlation between
the phases' random power imbalances.

The remainder of this paper is organized as follows. In Section \ref{sec:sys}, we describe the system model and formulate the optimization problem. In Section \ref{sec:id}, we propose both centralized and  distributed real-time algorithms for ideal energy storage. In Section \ref{sec:nonid}, we extend the algorithm to accommodate non-ideal energy storage with imperfect charging and discharging. Numerical results are presented in Section \ref{sec:sim}, and we conclude in Section \ref{sec:con}.

\section{System Model and Problem Statement}\label{sec:sys}
\begin{figure}[t]
\centering
\includegraphics[height= 1.3in,width = 3.2in]{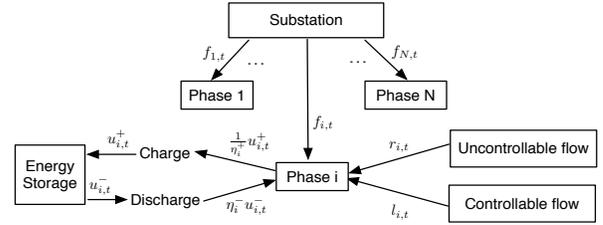}
\caption{System model with $N$ phases. The details of the $i$-th phase are shown.}
\label{fig:sys}
\end{figure}

Consider a discrete-time model with time  $t\in \{0,1,2, \ldots\}$. To simplify notation, we normalize the duration of each time period $\Delta t$ to one and thus eliminate  $\Delta t$ in presentation.
The system model is depicted in Fig. \ref{fig:sys}. A substation is connected with $ N\ge 2$ phases, each with single-phase loads and generation.\footnote{In practice, a typical substation consists of one or more three-phase feeders connected through a feeder breaker,  each of which  can supply multiple single-phase loads.  In this paper, since we focus on the problem of phase balancing in one feeder, the structure of feeders is omitted in Fig. 1. For a more thorough description of a distribution substation, please see Chapter 1 in \cite{bkdsma}.} We consider a general case where it is optional for each phase to deploy energy storage. Denote the set of  phases that deploy storage
by $\mathcal{E}\subseteq\{1, 2,\ldots,N\}$.
Below we first describe the components of each phase.

\subsection{System Model of Each Phase}\label{subsec:smph}
At the $i$-th phase, denote the amount of  uncontrollable power at time slot $t$ by $r_{i,t}$. The uncontrollable  flow can represent renewable generation such as wind and solar,  base loads, or the difference between renewable generation and base loads. Since the uncontrollable flow is generally governed by nature or uncertain human behavior, we assume that $r_{i,t}$ is random, but it is confined within an interval $[r_{i,\min}, r_{i,\max}]$.  Throughout the paper we use a bold letter to denote a vector that contains  elements  of $N$ phases. Here, we define $\br_t \define [r_{1,t},\ldots,r_{N,t}]$ to represent the uncontrollable flow vector at time slot $t$. The other vectors in the rest of this paper are defined similarly. 

Denote the amount of the controllable power flow at the $i$-th phase at time $t$ by $l_{i,t}$. The controllable flow can represent the output of conventional generators, or the consumption of flexible loads. We associate a cost function with the controllable flow and denote the function by $C_i(l_{i,t})$, which can represent the cost of local generators (e.g., an on-site diesel generator), or the cost of a utility for consuming power.

Denote the power flow between the substation and the $i$-th phase at time slot $t$  by $f_{i,t}$. 
Due to the capacity constraints of power lines,  the value of  $f_{i,t}$ is generally confined. We assume that at each time slot the power flow vector $\bff_t\in \mathcal{F}$, where the set $\mathcal{F}$ is non-empty, compact, and convex.
For example,  $\mathcal{F}$ may be defined as $\mathcal{F}\define\{\bff_t|f_{i,t}\in [f_{i,\min},f_{i,\max}], \forall i\}$.

\textbf{Remark: }The values of $r_{i,t}, l_{i,t},$ and $f_{i,t}$ can be positive or negative. 
We use the positive sign to indicate  power injection into the $i$-th phase, and the negative  sign to indicate  power extraction from the $i$-th phase. 

Assume that the $i$-th phase is equipped with an energy storage unit, i.e., $i\in\mathcal{E}$.
Denote the  charging and discharging rates of the storage at time slot $t$ by $u_{i,t}^+\in[0, u_{i,\max}]$ and $u_{i,t}^-\in[0, u_{i,\max}]$, respectively, where $u_{i,\max}$ is the maximum charging and discharging rates. Denote the energy state of the $i$-th storage at the beginning of time slot $t$ by $s_{i,t}$, which evolves as $s_{i,t+1} = s_{i,t} + u_{i,t}^+-u_{i,t}^-$.  The energy state $s_{i,t}$ is required to be within the storage's capacity limits $[s_{i,\min},s_{i,\max}]$. 

Due to conversion and storage losses,  charging and discharging  may not be perfectly efficient. For the $i$-th storage, we denote the charging efficiency  by $\eta_{i}^+\in (0,1]$ and the discharging efficiency by $\eta_{i}^-\in(0,1]$. Then, the associated charging and discharging quantities seen on each phase are $\frac{1}{\eta_i^+}u_{i,t}^+$ and $\eta_i^-u_{i,t}^-$, respectively (see Fig. \ref{fig:sys}). 
Owing to the round-trip efficiency or other operating constraints, 
simultaneous charging and discharging may be forbidden in practice, which can be reflected by the constraint $u_{i,t}^+ \cdot u_{i,t}^- =0, i\in\mathcal{E}$. 
Moreover, if the $i$-th phase is not equipped with storage, i.e.,  $i\notin\mathcal{E}$, we simply set the values of $s_{i,t}, u_{i,t}^+$, and $u_{i,t}^-$ to zero.

The energy storage can additionally be used for arbitrage.\footnote{Energy storage is still expensive based on the current technology. Therefore, in practice, besides phase balancing and energy arbitrage, energy storage can be used to provide other grid-wide services, e.g., volt-var control, load following, and peak shaving.}
Denote the electricity price at time slot $t$ by $p_{t}\in[p_{\min},p_{\max}]$, which is random over time. Then the cost of the $i$-th phase for energy arbitrage during time slot $t$ is $p_t(\frac{1}{\eta_i^+}u_{i,t}^+-\eta_i^-u_{i,t}^-)$.
Finally,  frequent charging and discharging can shorten the lifetime of storage\cite{ram02}. To model this effect, we introduce a degradation cost function $D_i(\cdot)$, with  negative input indicating discharging and positive input indicating charging.  Therefore, the degradation cost incurred at time slot $t$ is given by $D_{i}(u_{i,t}^+) +D_i(-u_{i,t}^-)$.\footnote{Accurate modeling of battery degradation is highly complicated and is  an active research area. In this paper, to focus on storage control,  we  employ a simplified degradation model, which is a function of  the charging and discharging amount.}

\subsection{Problem Statement}
Since phase imbalance is harmful for power system operation, it is critical to  balance the power flows $f_{i,t}$ among phases. To this end, we introduce a loss function $F(\cdot)$ to characterize the deviation of $f_{i,t}$ from the average power flow. In particular, for the $i$-th phase, $F(\cdot)$ is a function of  
$f_{i,t} - \overline{f}_{t}$, where $\overline{f}_{t}$ is the average defined as  $\overline{f}_{t} \define \frac{1}{N}\sum_{j=1}^Nf_{j,t}$. 

We assume that the system is operated by a representative of the substation, who aims to minimize the long-term system cost, which includes the costs of all phases. Specifically, based on the model described in Section \ref{subsec:smph}, the system cost at time slot $t$ is given by
\begin{align}
\nonumber
\textstyle{w_{t}} = & \textstyle{\sum_{i\in\mathcal{E}} \Big[p_t(\frac{1}{\eta_i^+}u_{i,t}^+-\eta_i^-u_{i,t}^-) + D_{i}(u_{i,t}^+) +D_i(-u_{i,t}^-)\Big]} \\
\nonumber
&\textstyle{+ \sum_{i=1}^N\big[C_i(l_{i,t})+F(f_{i,t} - \overline{f}_{t})\big].}
\end{align}

Denote the random system state at time slot $t$ by $\bq_t \define [\br_t,p_t]$, which includes the uncontrollable power flow of $N$ phases and the electricity price. Denote the control action at time slot $t$ by $\ba_t \define [\bl_t,\bu_t^+,\bu_t^-, \bff_t]$, which contains the controllable power flow, the charging and discharging amounts, and the power flow between each phase and the substation. We formulate the problem for phase balancing as the following stochastic optimization problem. 
\begin{align}
\nonumber
\textbf{P1:} \quad \min_{\{\ba_t\}}\quad &\limsup_{T\to\infty} \frac{1}{T}\sum_{t=0}^{T-1}\E[w_t]\\
\label{uit}
\st\quad & 0 \le u_{i,t}^+, u_{i,t}^-\le u_{i,\max},\forall i\in\mathcal{E},t,\\
\label{uitdc}
& u_{i,t}^+ \cdot u_{i,t}^- =0, \forall i\in\mathcal{E},t,\\
\label{sitt}
& s_{i,t+1} = s_{i,t} + u_{i,t}^+-u_{i,t}^-,\forall i\in\mathcal{E},t,\\
\label{sit}
& s_{i,\min}\le s_{i,t}\le s_{i,\max},\forall i\in\mathcal{E},t,\\
\label{uitn}
& u_{i,t}^- = u_{i,t}^+ =0, \forall i\notin\mathcal{E},t,\\
\label{fit}
& \bff_{t} \in \mathcal{F},  \forall t,\\
\label{bit}
&\hspace{-0.5cm} f_{i,t} + r_{i,t} +l_{i,t} + \eta_i^-u_{i,t}^- - \frac{1}{\eta_i^+}u_{i,t}^+= 0, \forall i, t.
\end{align}
The expectation on the objective is taken over the randomness of  $\bq_t$ and the possibly random control action that depends on $\bq_t$. Constraint \eqref{bit} enforces power balance at each phase at each time slot.

To keep mathematical exposition simple, we assume that the cost functions $C_i(\cdot)$ and $D_i(\cdot)$ are continuously differentiable and convex. This assumption is realistic because  many practical cost functions are well approximated this way\cite{bkpowerg}.  In particular, by convexity, the marginal cost is increasing. With the objective of minimizing the system cost, for the function $C(\cdot)$, this property discourages excessive use of the controllable flow.  For battery degradation, it is understood that faster charging or discharging has a more detrimental effect on the battery lifetime, and the convexity of $D(\cdot)$ reflects this behavior.
Denote  the derivatives of $C_i(\cdot)$ and $D_i(\cdot)$ by $C'_i(\cdot)$ and $D'_i(\cdot)$, respectively.  Since the variables $u_{i,t}^+, u_{i,t}^-,$ and $l_{i,t}$ are bounded based on the constraints of \textbf{P1}, the cost functions and their derivatives are bounded in the feasible set. For the cost function $C_i(\cdot)$, we denote its range by $[C_{i,\min},C_{i,\max}]$ and its range of the derivative by $[C'_{i,\min},C'_{i,\max}]$ in the feasible set. The range of the cost function $D_i(\cdot)$ and that of its derivative are defined similarly.
In addition, we assume that the loss function $F(\cdot)$ is convex and continuously differentiable.

We are interested in designing both centralized and distributed real-time algorithms for solving \textbf{P1}. Distributed implementation is motivated by the limited capability of real-time communication between the substation and each phase, and also the potential privacy concerns of each phase.
This is a challenging task due to   system uncertainty, the coupling of all phases through the objective and constraints, and the energy state constraint \eqref{sit} which couples the charging and discharging actions over time. In addition, we assume that the system is not equipped with  any forecaster and only has historical information of the system states. Designing appropriate forecasters and incorporating forecast into optimal control are important directions and are left for future work.

\section{Real-Time Algorithm for Ideal Energy Storage}\label{sec:id}
For tractability, in this section we first consider ideal energy storage   that has  perfectly efficient charging and discharging, i.e., $\eta_i^+ = \eta_i^- = 1$. The case of  non-ideal energy storage  is studied in Section \ref{sec:nonid}. 
We first propose  a centralized real-time algorithm that can be implemented by the substation and show its analytical performance. Then we provide distributed implementation for the proposed algorithm where  only limited information exchange is needed.

\subsection{Centralized Real-Time Algorithm and Analysis}
Under perfectly efficient charging and discharging, without loss of generality, we can combine the charging and discharging variables 
$u_{i,t}^+$ and $u_{i,t}^-$ into one by introducing a new variable $u_{i,t} \define u_{i,t}^+-u_{i,t}^-$,  which can represent the net charging and discharging amount. In particular, 
if $u_{i,t}>0$ it indicates charging, and if $u_{i,t}<0$ it indicates discharging. 

With the new  variable $u_{i,t}$, the non-simultaneous charging and discharging constraint \eqref{uitdc} can be eliminated, and the evolution of the energy state amounts to $s_{i,t+1} = s_{i,t} + u_{i,t}$. In addition, with $u_{i,t}$, the control action at time slot $t$ is now $\ba_t \define [\bl_t,\bu_t, \bff_t]$, and the system cost  can be rewritten as $w_{t} = \sum_{i\in\mathcal{E}} \big[p_tu_{i,t} + D_{i}(u_{i,t})\big] + \sum_{i=1}^N\big[C_i(l_{i,t})+F(f_{i,t} - \overline{f}_{t})\big]$.

For the design of real-time implementation, we employ Lyapunov optimization\cite{bkneely}, which has been used widely in wireless networks for dealing with time-averaged constraints and providing simple yet efficient algorithms for complex dynamic systems. However,  the energy state constraint \eqref{sit} is 
not a time-averaged constraint but a hard constraint, and it couples the control action $u_{i,t}$ over multiple time instances. As a result, \textbf{P1} is not amenable to the standard framework of Lyapunov optimization. 
To overcome this difficulty, 
we replace the energy state constraints \eqref{sitt} and \eqref{sit}  with a new time-averaged constraint, which only requires  the net charging and discharging amount to be zero on average, \ie,
\begin{align}\label{uitav}
\lim_{T\to\infty} \frac{1}{T}\sum_{t=0}^{T-1} \E[u_{i,t}] = 0, \forall i\in\mathcal{E}.
\end{align}
With the new constraint \eqref{uitav}, we form a  new  stochastic optimization problem as follows:
\begin{align}
	\nonumber
	\textbf{P2:} \quad \min_{\{\ba_t\}}\quad &\limsup_{T\to\infty} \frac{1}{T}\sum_{t=0}^{T-1}\E[w_{t}]\\
	\nonumber
	\st \quad &\eqref{fit},\;\eqref{uitav},\\
	&\label{biti} f_{i,t}+r_{i,t}+l_{i,t}-u_{i,t} = 0,\forall i,t,\\
	& \label{uitni} u_{i,t} = 0, \forall i\notin \mathcal{E}, t,\\
	&\label{uiti} -u_{i,\max}\le u_{i,t}\le u_{i,\max}, \forall i\in\mathcal{E},t.
\end{align}
It can be shown that  constraints \eqref{sitt} and \eqref{sit} imply \eqref{uitav}, i.e., any $u_{i,t}$ that satisfies \eqref{sitt} and \eqref{sit} also satisfies \eqref{uitav}. Hence, \textbf{P2} is a relaxation  of \textbf{P1} (see Appendix \ref{app:relax}). 

The above relaxation step is crucial for applying Lyapunov optimization. However, we emphasize that, solving \textbf{P2} is not our purpose (it is clear that, due to the relaxation of constraints \eqref{sitt} and \eqref{sit}, a solution to \textbf{P2} may be infeasible to \textbf{P1}). Instead, the significance of proposing \textbf{P2} is to facilitate the development of a real-time algorithm for \textbf{P1} and the associated performance analysis. Later we will prove in Proposition \ref{pro:beta} that our proposed algorithm  ensures that constraints \eqref{sitt} and \eqref{sit} are satisfied, and therefore produces a feasible solution to  \textbf{P1}.

We now propose a real-time algorithm leveraging Lyapunov optimization techniques. 
At  time slot $t$, for phase $i\in \mathcal{E}$, define a  Lyapunov function  
$L(s_{i,t}) \define  \frac{1}{2}(s_{i,t}-\beta_i)^2,$
which measures the deviation of the energy state $s_{i,t}$ from a perturbation parameter $\beta_i$.
The parameter $\beta_i$ is introduced to ensure the boundedness of the energy state, \ie, constraint \eqref{sit}, and it  needs to be carefully designed. In addition, 
we define a one-slot conditional Lyapunov drift as $\Delta(\bs_{t}) \define \E\big[\sum_{i\in\mathcal{E}}\frac{L(s_{i,t+1})-L(s_{i,t})}{V_i} |\bs_t\big]$, which collects a weighted  sum of the one-slot conditional drifts of the Lyapunov functions for all phases with storage. 

In our design of the real-time algorithm, instead of directly minimizing  the  system cost $w_t$, we consider a drift-plus-cost function  $\Delta(\bs_{t}) + \E[w_{t}|\bs_t]$. In particular, we first derive an upper bound on the drift-plus-cost function (see Appendix \ref{app:lemdpp} for the upper bound),
and then formulate a per-slot optimization problem to minimize this upper bound. Consequently, at each time slot $t$, we  solve the following optimization problem:
\begin{align*}
	\textbf{P3:} \quad \min_{\ba_t}\quad & w_t  + \sum_{i\in\mathcal{E}}\frac{(s_{i,t}-\beta_i)u_{i,t}}{V_i}\quad
	\st \; &\eqref{uitn}-\eqref{bit},\eqref{uiti}.
\end{align*}
Denote an optimal solution of \textbf{P3} at time slot $t$ by $\ba^*_t \define [\bl_t^*,\bu_t^{*}, \bff_t^*]$. At each time slot, after obtaining the solution $\ba^*_t$, we update $s_{i,t}$ using $u_{i,t}^*$.
It can be easily verified that the optimization problem \textbf{P3} is convex, and thus can be  solved by standard convex optimization software packages such as those in MATLAB. We will later shown in Theorem \ref{the:per} that such design of the per-slot optimization problem  leads to certain guaranteed performance.

In the  proposition below, we show that, despite the relaxation to arrive at $\textbf{P2}$, by appropriately designing the perturbation parameter $\beta_i$, we can ensure that constraint \eqref{sit} is satisfied, and therefore the control actions $\{\ba_t^*\}$ are feasible to \textbf{P1}.
\begin{pro}\label{pro:beta}
	For phase $i\in\mathcal{E}$, set the perturbation parameter $\beta_i$ as
	\begin{align}\label{betai}
		\beta_i\define s_{i,\min}+u_{i,\max}+V_i(p_{\max}+D'_{i,\max}+C'_{i,\max})
	\end{align}
	where  $V_i\in(0, V_{i,\max}]$ with 
	\begin{align}\label{vimax}
		\textstyle{V_{i,\max} \define \frac{s_{i,\max}-s_{i,\min}-2u_{i,\max}}{p_{\max}-p_{\min}+D'_{i,\max}-D'_{i,\min}+C'_{i,\max}-C'_{i,\min}}. }
	\end{align}
Then the control actions $\{\ba_t^*\}$ obtained by solving \textbf{P3} at each time $t$ are feasible to $\textbf{P1}$.
\end{pro}
\begin{IEEEproof}
	See Appendix \ref{app:probeta}.
\end{IEEEproof}

To ensure the positivity of $V_{i,\max}$ in \eqref{vimax}, we need  the numerator $s_{i,\max}-s_{i,\min}-2u_{i,\max}>0$. This is generally true for real-time applications, in which  the length of each time interval is small  ranging from a few seconds to minutes. 

The overall centralized real-time algorithm is summarized in Algorithm \ref{alg:rtid}, which can be implemented  by the substation.  It is worth mentioning that the proposed algorithm does not require any system statistics, which may be desirable when accurate system statistics are difficult to obtain.
\begin{algorithm}[t]
\caption{Centralized  algorithm for ideal storage.}
\label{alg:rtid}
At  time slot $t$, the substation executes the following steps sequentially:
\begin{algorithmic}[1]
\STATE observe the system state $\bq_t$ and the energy state $s_{i,t}$;
\STATE solve \textbf{P3} and obtain a solution $\ba^*_t \define [\bl_t^*,\bu_t^{*}, \bff_t^*]$; and
\STATE update $s_{i,t+1}$ by $s_{i,t} +u^*_{i,t}$. 
\end{algorithmic}
\end{algorithm}

Denote the optimal objective value of \textbf{P1} by $w^{\opt}$. Under Algorithm \ref{alg:rtid}, denote the objective value of \textbf{P1}  by $w^*$ and the system cost at time slot $t$ by $w_t^*$.
The performance of Algorithm \ref{alg:rtid} is stated in the following theorem.

\begin{theorem}\label{the:per}
Assume that the system state $\bq_t$ is i.i.d. over time and the equipped storage at the phases is perfectly efficient. 
Under Algorithm \ref{alg:rtid} the following statements hold.
\begin{enumerate}
\item $w^*-w^{\opt}\le   \sum_{i\in\mathcal{E}}\frac{u_{i,\max}^2}{2V_i}$.
\item $\frac{1}{T}\sum_{t=0}^{T-1}\E[{w}_t^*] - w^{\opt}\le  \sum_{i\in\mathcal{E}}\frac{u_{i,\max}^2}{2V_i} + \frac{\E[L(s_{i,0})]}{TV_i}.$
\end{enumerate}
\end{theorem}
\begin{IEEEproof}
See Appendix \ref{app:theper}.
\end{IEEEproof}
\textbf{Remarks:} 
\begin{itemize}
\item For Theorem \ref{the:per}.1, first, if $\mathcal{E}$ is empty, \ie, no phase deploys storage, then Algorithm \ref{alg:rtid} achieves the optimal objective value. In fact, for this case, Algorithm \ref{alg:rtid} reduces to a greedy algorithm that only minimizes the current system cost  at each time. Second, if $\mathcal{E}$ is non-empty, to minimize the gap to the optimal objective value, we should set $V_i = V_{i,\max}$. Asymptotically, if the energy capacity $s_{i,\max}$ is large and thus $V_{i,\max}$ is large,  Algorithm \ref{alg:rtid} achieves the optimal objective value.\footnote{The choice of $V_i = V_{i,\max}$ and the asymptotic optimality are based on the linear storage model in \eqref{sitt}. These conclusions need to be re-examined when a more general storage model with other factors such as storage efficiency is considered \cite{qcyr14,qcyr14b}.}
\item 
In Theorem \ref{the:per}.2, we characterize the performance of Algorithm \ref{alg:rtid} over a finite time horizon.   The result not only shows the performance gap of the algorithm over a finite time $T$, but also reveals how the gap converges asymptotically to the one in Theorem \ref{the:per}.1 as $T$ grows. It can be seen that, the gap contains a component $\sum_{i\in\mathcal{E}}\frac{\E[L(s_{i,0})]}{TV_i}$  due to the initialization of the energy states, which linearly decreases with the time horizon $T$.
\item The i.i.d. assumption of the system state $\bq_t$ can be relaxed to accommodate $\bq_t$ that  follows a finite state irreducible and aperiodic Markov chain. Using a multi-slot drift technique\cite{bkneely}, we can show similar conclusions which are omitted here. In simulation, we will evaluate the algorithm performance when the uncontrollable power flows are temporally correlated. 
\end{itemize}

An interesting additional consequence of Theorem \ref{the:per} is that we obtain a general rule of thumb for the allocation of energy storage capacity among the phases.
In particular, in the following proposition, we demonstrate that, under some mild conditions, equal allocation of a given energy storage capacity  results in a lower overall system cost.
\begin{pro}\label{pro:simax}
Assume that $s_{i,\min}, u_{i,\min}$ and $D_{i,\max}-D_{i,\min}$ are identical for all $i\in \mathcal{E}$. Assume further that for all phases, $C_{i,\max}'-C_{i,\min}'$ is the same. Then, under Algorithm \ref{alg:rtid}, if the total energy storage capacity $\sum_{i\in \mathcal{E}} s_{i,\max}$ is fixed and the control parameter $V_i = V_{i,\max}$ is as in \eqref{vimax}, the upper bound of the performance gap in Theorem \ref{the:per}.1, i.e., $\sum_{i\in\mathcal{E}}\frac{u_{i,\max}^2}{2V_i}$, is minimized when the energy storage capacity is equally allocated among phases.
\end{pro}
\begin{IEEEproof}
See Appendix \ref{app:prosimax}.	
\end{IEEEproof}

The above result states that energy storage is best allocated equally over the phases. Note that this result is robust because it does not depend on any system statistics or specific values of  system parameters. We will revisit this  in simulation.

In this paper, as our focus is on   designing real-time algorithms for storage control,  we  use a stylized  system model shown in Fig. \ref{fig:sys}. In particular, we do not model the  network structure at each phase, and therefore do not consider how to place storage. Instead, we assume a given arbitrary deployment of storage at any location of each phase. Nonetheless, we point out that storage placement  affects the investment strategy of  power systems and is crucial for grid operation. This problem has attracted considerable attention and has been investigated in many papers (e.g., \cite{tbh13}). In particular, in the case of  physical storage, the authors of \cite{tbh13} proved that, under some technical conditions, there always exists an optimal strategy of storage placement that assigns zero storage at generation-only buses that link to the rest of the network via single transmission lines. Moreover, if storage is  load aggregation, then it can be distributed over the phase. How to extend our results to consider   storage placement is a topic for  future study.

\subsection{Distributed  Implementation of Centralized  Algorithm}\label{subsec:dis}
To accomplish the  implementation of Algorithm \ref{alg:rtid} in a centralized way, each phase has to provide all information that is required to solve the real-time problem \textbf{P3}. Specifically,  for each phase, 
the cost functions and the associated optimization constraints need to be communicated to the substation in advance. In addition, at each time slot,  the information of the uncontrollable power flow as well as the storage energy state has to be sent to the substation. 
However, in practice, due to the limited capability of real-time communication along with  potential privacy concerns of each phase, some of the  aforementioned information may be unavailable at the substation. Therefore, the centralized implementation may be infeasible.  In this subsection, we provide a distributed algorithm for solving \textbf{P3} in which only limited information exchange  is required.
For ease of notation,  we suppress the time index $t$ in the following presentation.

The distributed algorithm is based on  the alternating direction method of multipliers (ADMM)  \cite{adboyd}. 
To facilitate algorithm development, we rewrite \textbf{P3} as follows:
\begin{align}
\nonumber
\min_{\ba}\quad &\mb{1}(\bff\in \mathcal{F}) + \sum_{i=1}^N\Big[H_i(l_i,u_i)  +F(f_i-\overline{f})\Big]\\
\label{bi}
\st \quad & f_{i} + r_{i} +l_{i} -u_{i}= 0, \forall i
\end{align}
where $\mb{1}(\cdot)$ is the indicator function that equals $0$ (resp. $+\infty$) when the enclosed event is true (resp. false), and for each phase the function $H_{i}(l_i,u_i)$ is defined as follows:
\begin{align*}
H_{i}(l_i,u_i) \define
\begin{cases}
\frac{(s_{i}-\beta_i)u_i}{V_i}+ pu_i+D_i(u_i)+ C_i(l_{i}) \\
+ \mb{1}(-u_{i,\max}\le u_{i}\le u_{i,\max}), & \textrm{ if } i\in\mathcal{E}\\
 C_i(l_{i}) + \mb{1}(u_i = 0), & \textrm{ if } i\notin\mathcal{E}.
\end{cases}
\end{align*}
We associate a Lagrange multiplier $\lambda_i$ with equality \eqref{bi}.

By treating the variables $(\bl,\bu)$ as one block and the variable $\bff$ as the other, we  express
the updates at the $(k+1)$-th iteration below according to the ADMM algorithm.

\vspace{-0.2cm}
{\small{
\begin{align*}
&(l_i,u_i)^{k+1}\leftarrow  \argmin_{l_i,u_i}\Big[H_i(l_i,u_i)+\frac{\rho}{2}(f_{i}^k + r_{i} +l_{i}  -u_{i}+\frac{\lambda_i^k}{\rho})^2 \Big]\\
&\bff^{k+1}\leftarrow \argmin_{\bff\in\mathcal{F}}\sum_{i=1}^N\Big[F(f_i-\overline{f})+\frac{\rho}{2}(f_{i} + r_{i} +l_{i}^{k+1}  -u_{i}^{k+1}+\frac{\lambda_i^k}{\rho})^2\Big]\\
&\lambda_i^{k+1}\leftarrow \lambda_i^k + \rho(f_{i}^{k+1} + r_{i} +l_{i}^{k+1}  -u_{i}^{k+1})
\end{align*}}}
\hspace{-0.2cm}where $\rho>0$ is a pre-determined parameter.

\begin{figure}[t]
\centering
\includegraphics[height = 1in, width = 2.2in]{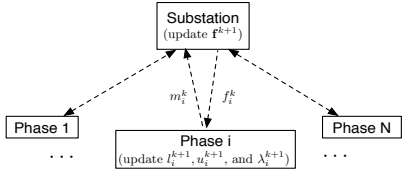}
\caption{Distributed implementation for solving \textbf{P3}.}
\label{fig:dis}
\end{figure}

To implement the above iteration, each phase  updates the controllable power flow $l_i$, the net charging and discharging amount $u_i,$ and the Lagrange multiplier $\lambda_i$, while the substation updates the power flow vector $\bff$. 
For  information exchange, at the $(k+1)$-th iteration, the substation sends $f_{i}^k$ to each phase, and each phase  provides $m_i^k \define r_{i} +l_{i}^{k+1} -u_{i}^{k+1}+\frac{\lambda_i^k}{\rho}$ to the substation. The schematic representation of the distributed implementation  is given in Fig. \ref{fig:dis}.

\textbf{Remark:} Although the communication network structure is the same in both centralized and distributed algorithms (i.e., star topology),  with the proposed distributed algorithm, each phase only needs to provide the update of $m_i^k$ to the substation without revealing the cost functions or  other parameters. Therefore, the communication load and the information revealed by each phase are limited.

The convergence behavior of the distributed algorithm is summarized in the following theorem. The proof follows Theorem 2 in \cite{wb13} and thus is omitted.

\begin{theorem}
Assume that the functions $D_i(\cdot), C_i(\cdot),$ and $F(\cdot)$ are closed, proper, and convex.  The sequence $\{\bl^k, \bu^k, \bff^k, \lambda^k\}$ converges to an optimal primal-dual solution of \textbf{P3} with the worst case convergence rate $O(1/k)$.
\end{theorem}

\section{Extension to Non-ideal Energy Storage}\label{sec:nonid}
In this section, we  discuss the algorithm design for  non-ideal energy storage with inefficient charging and discharging. This is significant because common storage technologies such as batteries can have  round-trip efficiency, i.e., $\eta_i^+\cdot\eta_i^-$, ranging from $70\%$ to $95\%$\cite{cg14}.

The mathematical framework of the algorithm design follows that of ideal storage. However, due to imperfect charging and discharging, the charging and discharging variables $u_{i,t}^+$ and $u_{i,t}^-$ cannot be combined into one as we did in Section \ref{sec:id}, and therefore, the (non-convex) non-simultaneous charging and discharging constraint \eqref{uitdc} cannot be eliminated. 
To overcome this difficulty, we first ignore constraint \eqref{uitdc} and then adjust the resultant solution to satisfy the constraint.  

Specifically, we first modify the per-slot optimization problem \textbf{P3}  to the following:
\begin{align}
\nonumber
\textbf{P3':} \quad \min_{\ba_t}\quad & w_t  + \sum_{i\in\mathcal{E}}\frac{(s_{i,t}-\beta_i)}{V_i}(u_{i,t}^+-u_{i,t}^-)\\
\nonumber
\st \quad &\eqref{uit},\eqref{uitn}-\eqref{bit}
\end{align}
where we have defined the perturbation parameter 
\begin{align}\label{betain}
	\textstyle{\beta_i\define s_{i,\min}+u_{i,\max}+V_i(\frac{p_{\max}}{\eta_i^+}+\frac{1}{\eta_i^+}C'_{i,\max}+D'_{i,\max}).}
\end{align}
The parameter  $V_i$ in \eqref{betain} lies in the interval $(0, V_{i,\max}]$, where 
	$V_{i,\max} \define \frac{s_{i,\max}-s_{i,\min}-2u_{i,\max}}{\frac{p_{\max}}{\eta_i^+}-p_{\min}\eta_i^-+D'_{i,\max}-D'_{i,\min}+\frac{1}{\eta_i^+}C'_{i,\max}-\eta_i^-C'_{i,\min}}. $
Note that the definition of  $\beta_i$ in \eqref{betain} is similar to that in \eqref{betai} for ideal storage, except the inclusion of the charging and discharging efficiencies.  Moreover, if $\eta_i^+ =\eta_i^- = 1$,  \eqref{betain} reduces to \eqref{betai}.

The overall centralized algorithm is summarized in Algorithm \ref{alg:rtnid}, where we use the superscript notations $\hat{}$ and $*$ to indicate the intermediate solution derived from \textbf{P3'} and the final solution, respectively.
To ensure that the final solution satisfies constraint \eqref{uitdc},   
in Step 3, we adjust the  intermediate charging and discharging solutions  $\hat{u}_{i,t}^+$ and $\hat{u}_{i,t}^{-}$, and  the controllable power flow $\hat{l}_{i,t}$, so that simultaneous charging and discharging cannot happen and the power balance constraint \eqref{bit} still holds. 

\textbf{Remarks:} Under some conditions, constraint \eqref{uitdc} may automatically hold by solving \textbf{P3'}, e.g., when the electricity price $p_t$ is positive and the cost function of the controllable flow $C_i(\cdot)$ is increasing. However, if $p_t$ can be negative or consuming controllable flow costs money, the solution of  \textbf{P3'} may not meet  constraint  \eqref{uitdc}  and thus Step 3 in Algorithm \ref{alg:rtnid} may be necessary. In addition, if simultaneous charging and discharging is allowed in practice,  we can simply eliminate Step 3 in Algorithm \ref{alg:rtnid}.

\begin{algorithm}[t]
\caption{Centralized  algorithm for non-ideal storage.}
\label{alg:rtnid}
At  time slot  $t$, the substation executes the following steps sequentially:
\begin{algorithmic}[1]
\STATE observe the system state $\bq_t$ and the energy state $s_{i,t}$;
\STATE solve \textbf{P3'} and obtain an intermediate solution $\hat{\ba}_t \define [\hat{\bl}_t,\hat{\bu}_t^+,\hat{\bu}_t^-, \hat{\bff}_t]$;
\STATE generate the final solution ${\ba}_t^*$ where $u_{i,t}^{+*} = \max\{\hat{u}_{i,t}^+-\hat{u}_{i,t}^-,0 \}, u_{i,t}^{-*} = \max\{\hat{u}_{i,t}^--\hat{u}_{i,t}^+,0 \}, $ $l_{i,t}^* = \hat{l}_{i,t} + \eta_i^-\hat{u}_{i,t}^- - \frac{1}{\eta_i^+}\hat{u}_{i,t}^+ - \eta_i^-{u}_{i,t}^{-*} + \frac{1}{\eta_i^+}{u}_{i,t}^{+*}$, and $\bff^*_t = \hat{\bff}_t$; and
\STATE update $s_{i,t}$ by \eqref{sitt} using $u_{i,t}^{+*}$ and $u_{i,t}^{-*}$.
\end{algorithmic}
\end{algorithm}

The performance of  Algorithm \ref{alg:rtnid} is summarized in the following theorem.
\begin{theorem}\label{the:nonper}
Assume that the system state $\bq_t$ is i.i.d. over time and the equipped storage at the phases is not perfectly efficient. 
Under Algorithm \ref{alg:rtnid} the following statements hold.
\begin{enumerate}
\item $\{\ba_t^*\}$ is feasible for \textbf{P1}.
\item $w^*-w^{\opt}\le   \sum_{i\in\mathcal{E}}\frac{u_{i,\max}^2}{2V_i}+\epsilon$.
\item $\frac{1}{T}\sum_{t=0}^{T-1}\E[{w}_t^*] - w^{\opt}\le   \sum_{i\in\mathcal{E}}\Big[\frac{u_{i,\max}^2}{2V_i} + \frac{\E[L(s_{i,0})]}{TV_i} \Big]+\epsilon,$
\end{enumerate}
where $\epsilon\define \sum_{i\in\mathcal{E}}p_{\max}u_{i,\max}(\frac{1}{\eta_i^+}+\eta_i^-)+2D_{i,\max}+C_{i,\max}$.
\end{theorem}
\begin{IEEEproof}
See Appendix \ref{app:thenonper}.
\end{IEEEproof}

The results in Theorem \ref{the:nonper} parallel those in Theorem \ref{the:per} for ideal storage, with an extra gap $\epsilon$ incurred due to the adjustment of the intermediate solutions. Furthermore, since constraint \eqref{uitdc} is ignored in  \textbf{P3'}, the problem is  convex and therefore Algorithm \ref{alg:rtnid} can be implemented distributively  using a similar  ADMM-based algorithm as that in Section \ref{subsec:dis}.

\section{Numerical Results}\label{sec:sim}
In this section, we numerically evaluate the performance of the proposed algorithm. In each example, all  phases  are equipped with energy storage.
The specific values for the system parameters and functions are shown in Table \ref{tab:setup}.
The other default setup is as follows: the system state $[\br_t,q_t]$ is i.i.d. over time; at each time slot, the uncontrollable power flows are modeled as independent among phases, and they follow the Gaussian distribution $\mathcal{N}(0, 4^2)$ truncated within $[r_{i,\min},r_{i,\max}]$; and  the electricity price $p_t$ is approximated to follow the uniform distribution. 
For ideal storage, at each time slot, the control action $\ba_t \define [\bl_t, \bu_t, \bff_t]$ is generated by Algorithm \ref{alg:rtid}, and for non-ideal storage, $\ba_t\define [\bl_t, \bu_t^+,\bu_t^-, \bff_t]$ is generated by Algorithm \ref{alg:rtnid}.
Both Algorithms  are run for $T = 500$ time slots. The control parameter $V_i$ is set to $V_{i,\max}$.  Note that the only difference between the centralized and the distributed algorithms is whether the per-slot optimization problem in Algorithm \ref{alg:rtid} (or Algorithm \ref{alg:rtnid}) is solved centrally by the substation or  in distributed fashion by the substation and all phases. Therefore, both algorithms lead to the same solution to the per-slot optimization problem and thus the same  time-averaged system cost.\footnote{Since our focus in this paper is the design of real-time algorithms for storage control, we  use a stylized system model.  Extending our algorithm to accommodate more details of a power system and  implementing the algorithm in a real network using full transient simulation are topics of future work.}

For comparison, we use a greedy algorithm as the benchmark, which does not account for  the future performance. In particular, at each time slot, the greedy algorithm  minimizes the current system cost  subject to all constraints of \textbf{P1}. For ideal storage, the greedy algorithm solves the following optimization problem in each time slot:
\begin{align*}
	&\min_{\bl_t,\bu_t,\bff_t}\quad w_t\\
	&\st \quad \eqref{fit}, \eqref{biti}, \eqref{uitni}, \\
	&\hspace{.8cm} u_{i,t}\ge \max\{-u_{i,\max},s_{i,\min}-s_{i,t}\}\\
	& \hspace{.8cm}  u_{i,t}\le \min\{u_{i,\max},s_{i,\max}-s_{i,t}\}.
\end{align*}
For non-ideal storage, at time slot $t$, an intermediate solution is first found by solving the optimization problem 
\begin{align*}
	&\min_{\bl_t,\bu_t^+,\bu_t^-,\bff_t}\quad w_t \quad \quad \st \quad \eqref{uit},\eqref{sitt}-\eqref{bit}
\end{align*}
without the non-simultaneous charging and discharging constraint \eqref{uitdc}. Then, the final solution of the greedy algorithm is determined by adjusting the intermediate solution using Step 3 in Algorithm \ref{alg:rtnid}.

\begin{table}[t]
\renewcommand{\arraystretch}{1}
\caption{Default setup of parameters and functions}
\label{tab:setup}
\centering
\begin{tabular}{p{1.8cm}p{2.5cm}|p{1.2cm}p{0.8cm}}
\hline
Par.  & Setup & Par. (Fun.) & Setup  \\
\hline
\hline
$[r_{i,\min}, r_{i,\max}]$ & $[-8,8]$ (kW) & $\eta_{i}^+, \eta_i^-$ & $1$\\
$[f_{i,\min},f_{i,\max}]$ & $[-5,5]$ (kW)  &  $C_i(x)$ & $1.5x^2$\\
$[s_{i,\min},s_{i,\max}]$ & $[2, 10]$ (kWh) &$D_i(x)$ & $0.2x^2$\\
$[p_{\min},p_{\max}]$ & $[7,12]$ (cents/kWh) & $F(x)$ & $10x^2$\\
$u_{i,\max}$  & $1$ (kW) & N& 3\\
\hline
\end{tabular}
\end{table}

\subsection{Effect of  Energy Capacity of Storage}
\begin{figure}[t]
\centering
\includegraphics[height= 2in,width = 2.5in]{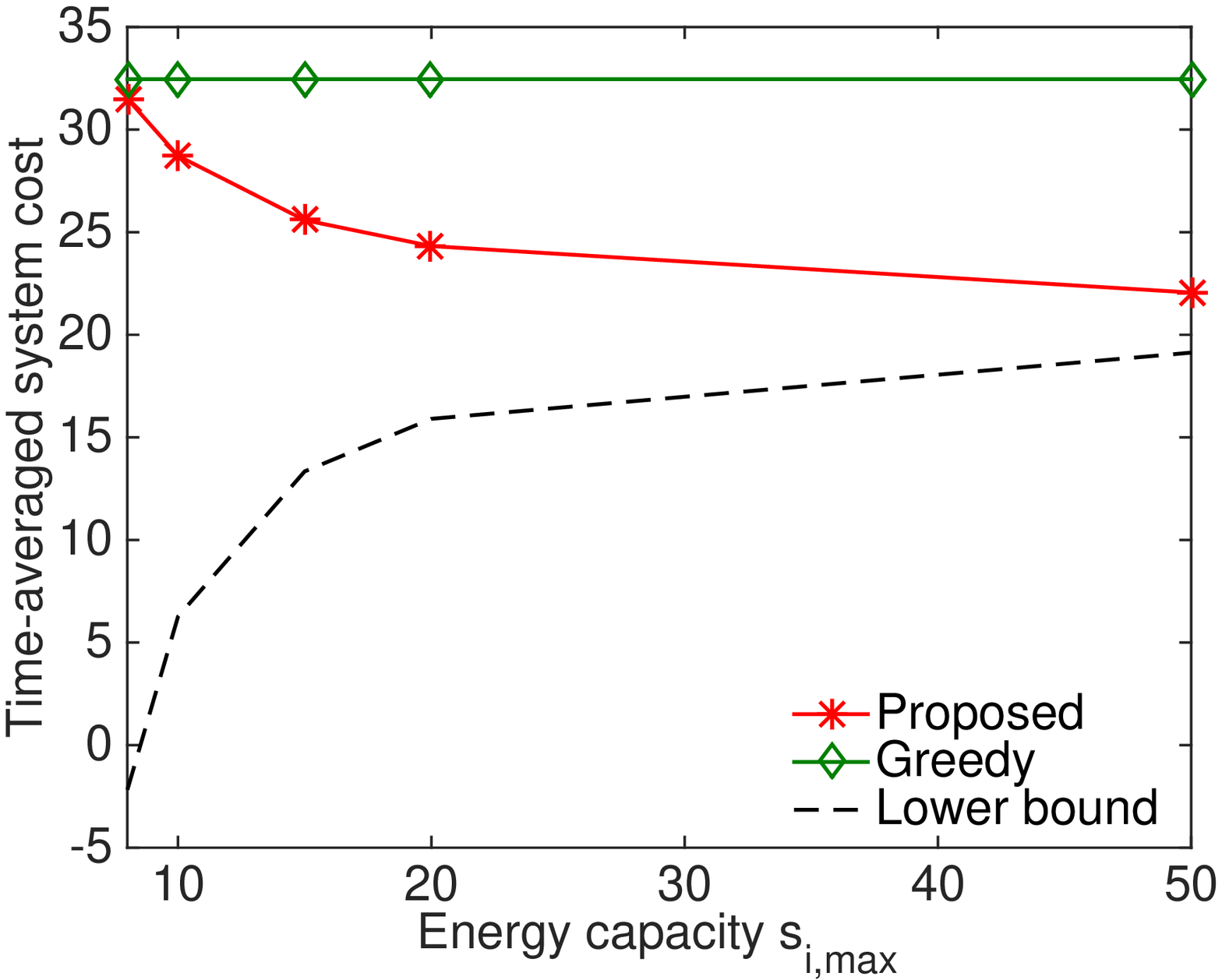}
\caption{System cost versus energy capacity of storage ($s_{1,\max} = s_{2,\max} = s_{3,\max}$).}
\label{fig:difsmax}
\centering
\includegraphics[height= 2in,width = 2.5in]{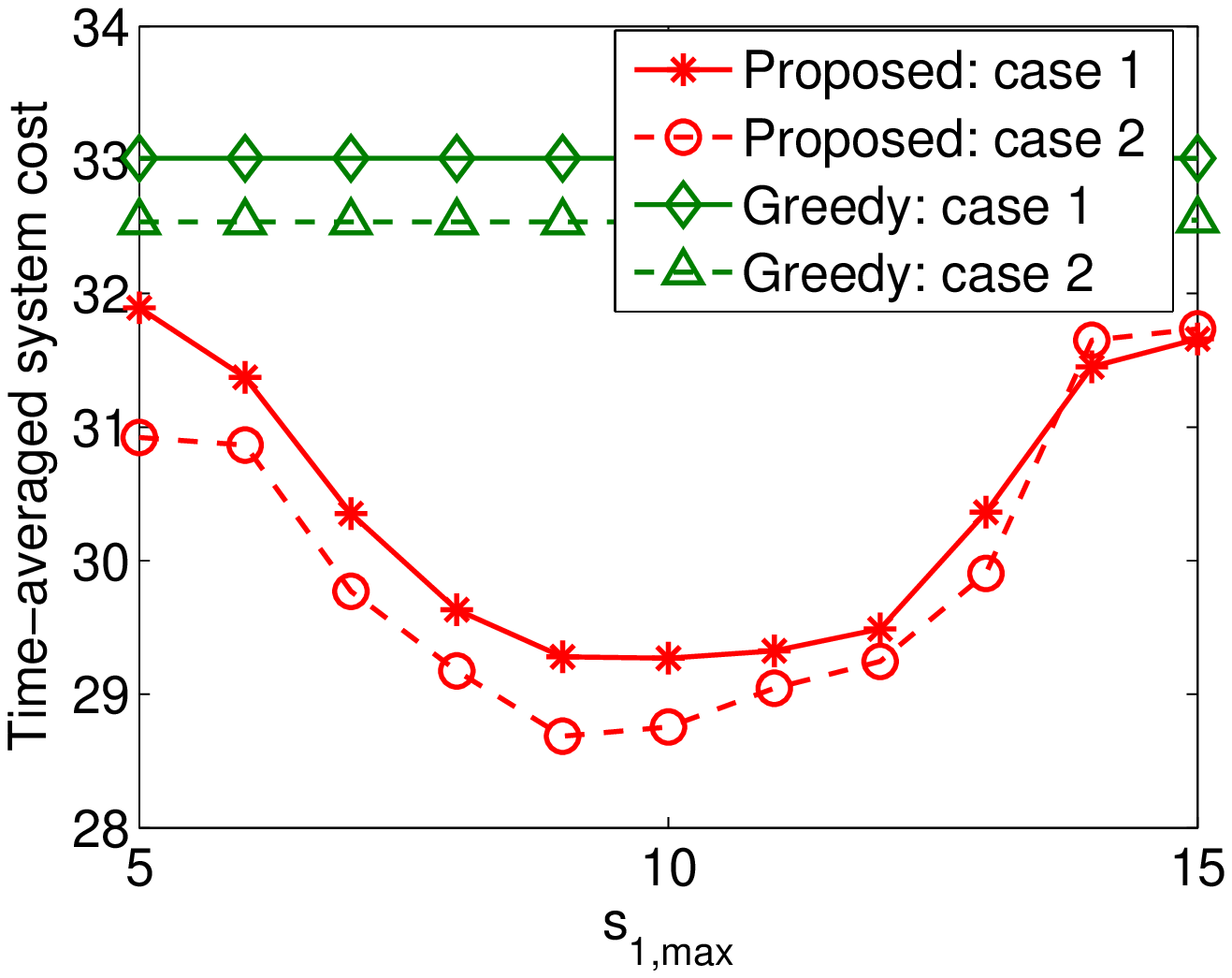}
\caption{System cost versus energy capacity of storage at Phase 1  ($s_{1,\max} + s_{2,\max} + s_{3,\max} = 30$ kWh).}
\label{fig:compdifsmax}
\end{figure}
	
In this subsection, we consider the effect of  energy capacity allocation on the system cost.
In Fig. \ref{fig:difsmax}, we increase the values of the energy capacity of all storage units from $8$ kWh to $50$ kWh. Note that for the proposed algorithm the role of $s_{i,\max}$ is played through the design of the control parameter $V_{i,\max}$ in \eqref{vimax}, and for the greedy algorithm the effect of $s_{i,\max}$ is reflected through the upper bound of the net charging and discharging variable $u_{i,t}$ in the optimization problem. We see that, as $s_{i,\max}$ increases, the system cost of the greedy algorithm does not change, while that of the proposed algorithm drops with a decreasing slope. The former phenomenon could happen when the maximum charging and discharging rate  $u_{i,\max}$ is relatively small and thus $s_{i,\max}$ has limited effect on $u_{i,t}$.
The latter observation is consistent with the second remark below Theorem \ref{the:per} that the proposed algorithm is asymptotically optimal when $s_{i,\max}$ is large.   In addition, from Theorem \ref{the:per}.1 we can obtain a lower bound of the minimum system cost as $w^{\textrm{opt}} \ge w^*-\sum_{i\in\mathcal{E}}\frac{u_{i,\max}^2}{2V_i}$. In Fig. \ref{fig:difsmax}, we also show the curve of this lower bound. In particular, when the energy capacity is large, this lower bound is tight. However, when the energy capacity is small-to-moderate, this lower bound is loose. For the remaining part of the simulation section, we study the performance of the proposed algorithm when the energy capacity is moderate (e.g., $s_{i,\max} = 10$ kWh).  Therefore, we have omitted this low bound in the remaining figures. Instead, the benchmark greedy algorithm is more effective in evaluating the numerical performance of the proposed algorithm. Moreover, like the proposed algorithm, the performance of the greedy algorithm  serves as an upper bound of the minimum system cost.

In Fig. \ref{fig:compdifsmax}, we fix the total energy capacity of all storage units to $30$ kWh (i.e., $s_{1,\max} + s_{2,\max} + s_{3,\max} = 30$ kWh) and vary the capacity allocation among phases. In particular, we fix $s_{2,\max}$ at $10$ kWh and change $s_{1,\max}$  from $5$ kWh to $15$ kWh.   Two cases are considered: Case 1, the standard deviation of the uncontrollable flow of each phase is $4$ kW (default setup); Case 2, the standard deviations of the uncontrollable flow of phases 1, 2, and 3 are $3$ kW, $4$ kW, and $5$ kW, respectively. For both algorithms, Case 2 leads to a smaller system cost in general. Moreover, for the greedy algorithm, the system cost barely changes with  $s_{1,\max}$.  In comparison,  for the proposed algorithm, the system cost achieves the lowest value when the energy capacity is approximately equally allocated. This observation is consistent with our conclusion in Proposition \ref{pro:simax}. 

\subsection{Effect of Correlations of Uncontrollable Power Flows}
\begin{figure}[t]
\centering
\includegraphics[height= 2in,width = 2.5in]{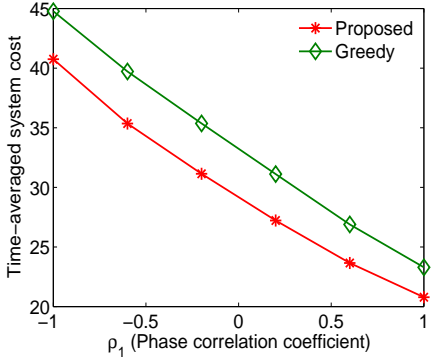}
\caption{System cost versus phase correlation coefficient of uncontrollable power flows of Phase 1 and Phase 2.}
\label{fig:phase_cor1}
\end{figure}
\begin{figure}[t]
\centering
\includegraphics[height= 2in,width = 2.5in]{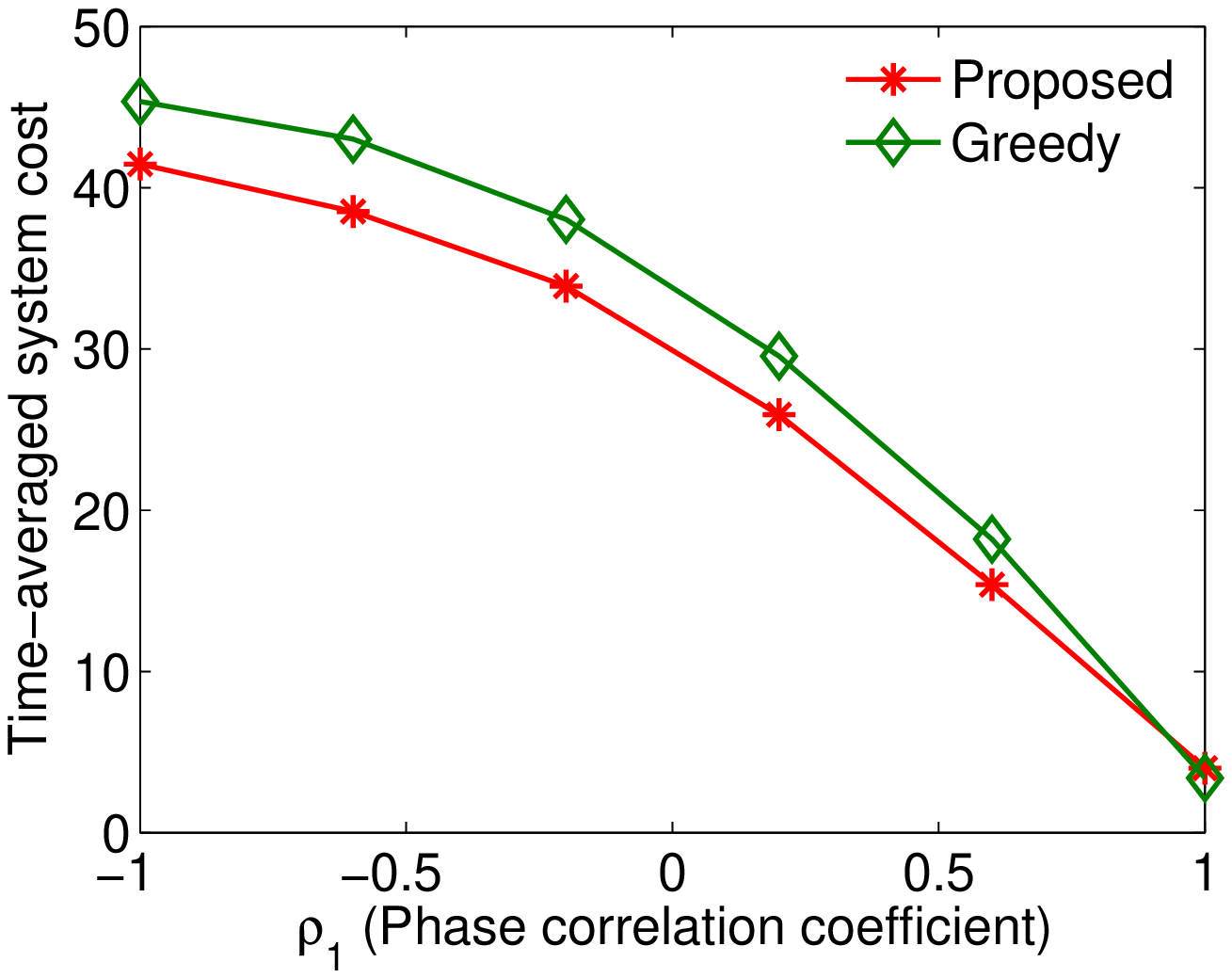}
\caption{System cost versus phase correlation coefficient of uncontrollable power flows of Phase 1 and Phases 2, 3.}
\label{fig:phase_cor2}
\end{figure}
\begin{figure}[t]
\centering
\includegraphics[height= 2in,width = 2.5in]{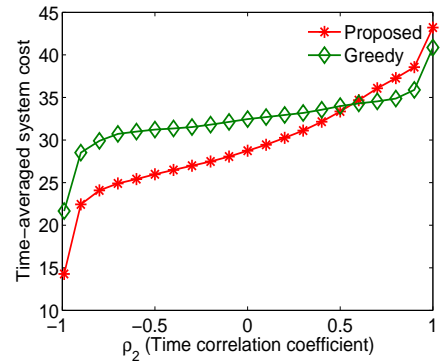}
\caption{System cost versus time correlation coefficient of uncontrollable power flow at each phase.}
\label{fig:markovflow}
\end{figure}
\begin{figure}[t]
\centering
\includegraphics[height= 2in,width = 2.5in]{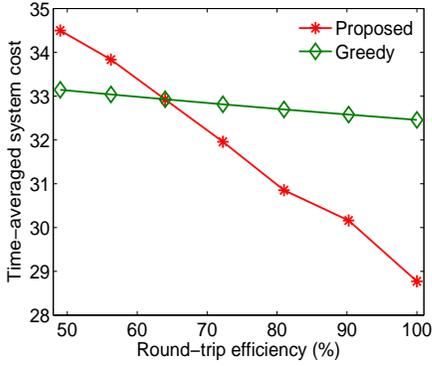}
\caption{System cost versus round-trip efficiency of storage at each phase.}
\label{fig:difcdcoef}
\end{figure}

In this subsection, we examine the effect of both the phase and time correlation of the uncontrollable power flows on the system cost. In Fig. \ref{fig:phase_cor1}, we assume that at each time slot, the uncontrollable flows of Phases 1 and 2 are correlated with the phase correlation coefficient, denoted by $\rho_1$, while the uncontrollable flow of Phase 3 is independent of those of Phases 1 and 2. We  see that, for both algorithms, the system cost decreases with $\rho_1$. This is easy to understand, since with a larger $\rho_1$ the uncontrollable flows of Phases 1 and 2 are more positively related, which makes phase balancing less challenging. In Fig. \ref{fig:phase_cor2}, we additionally assume that the uncontrollable flow of Phase 3 is correlated with that of Phase 1 with the same correlation coefficient $\rho_1$. With the additional correlation among phases, the performance gap   between the proposed algorithm and the greedy algorithm becomes smaller.  

In Fig. \ref{fig:markovflow}, we assume that the uncontrollable flows are independent among phases at each time slot, but they are  temporally correlated with the time correlation coefficient, denoted by $\rho_2$. We observe that, for both algorithms, the system cost  increases with $\rho_2$. This is because  at each phase, when the uncontrollable flow is more positively correlated, the more expensive controllable flow is used for phase balancing since the energy state of the storage is close to its range limit. Consequently, the proposed algorithm achieves a lower system cost when the uncontrollable flow is more negatively correlated.

\begin{figure}[t]
\centering
\includegraphics[height= 2in,width = 2.5in]{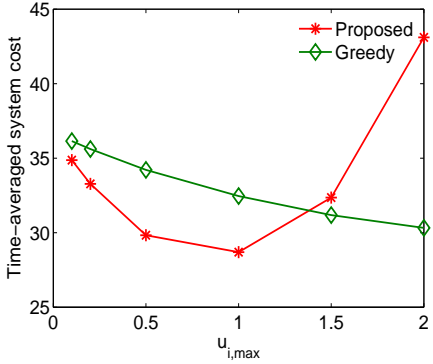}
\caption{System cost versus maximum charging/discharging rate of storage at each phase.}
\label{fig:difumax}
\end{figure}

\begin{figure}[t]
\centering
\includegraphics[height= 2in,width = 3.5in]{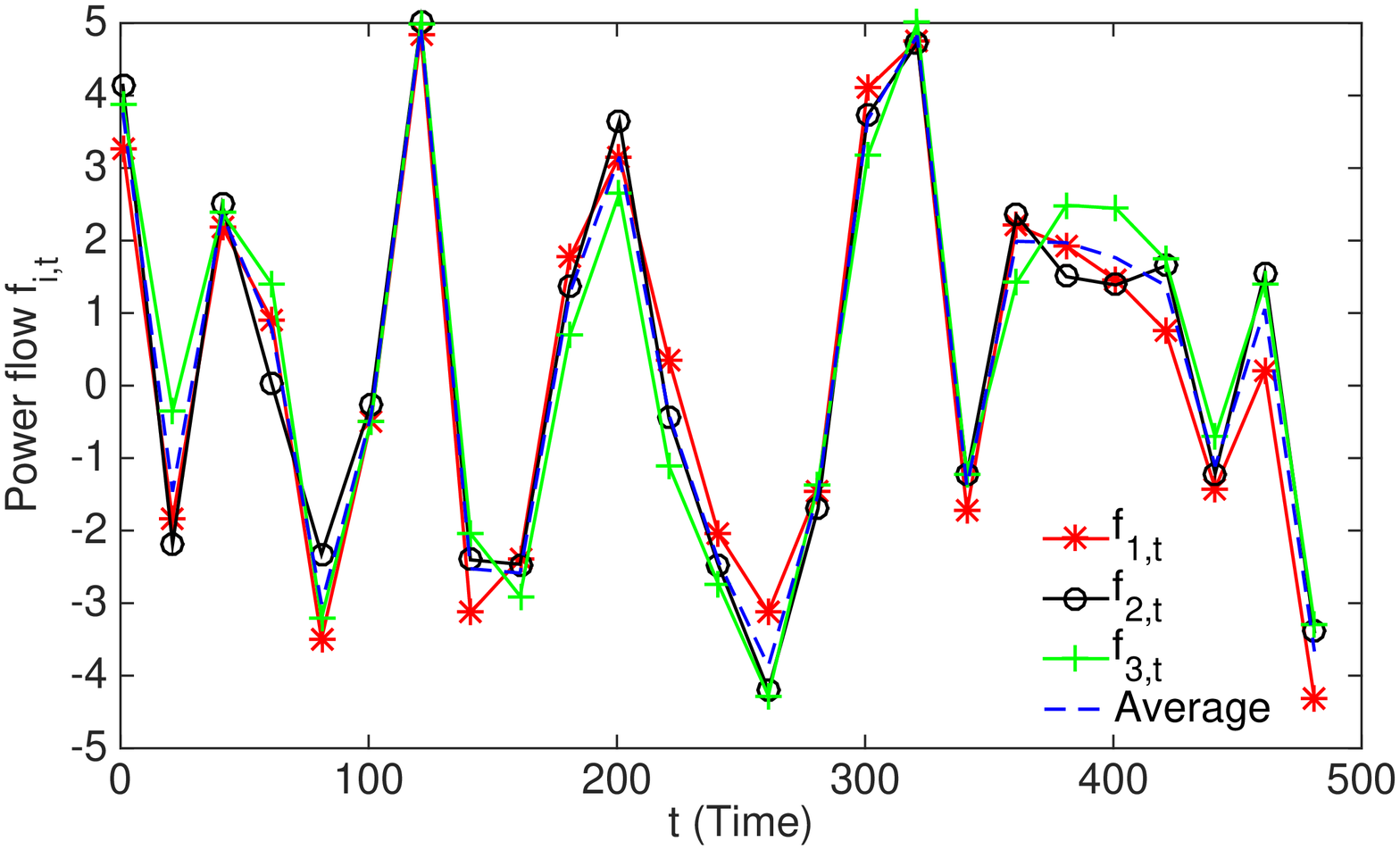}
\caption{Time trajectory of power flow $f_{i,t}$, for $i=1,2,3$ .}
\label{fig:power_flow}
\end{figure}
\begin{figure}[t]
\centering
\includegraphics[height= 2in,width = 2.5in]{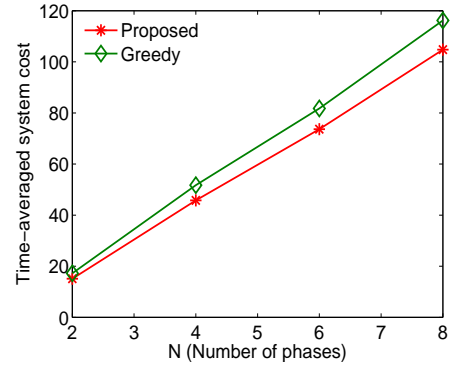}
\caption{System cost versus number of total phases. }
\label{fig:difn}
\end{figure}

\subsection{Effect of Charging and Discharging Circuit Parameters}
In Fig. \ref{fig:difcdcoef}, we consider that each phase is equipped with non-ideal energy storage. The charging and discharging efficiencies $\eta_i^+$ and $\eta_i^-$ of each storage are assumed to be the same. We see that for both algorithms, the system cost decreases almost linearly with the round-trip efficiency. The decreasing trend is expected since the storage becomes more efficient with a larger value of the round-trip efficiency. In particular,  the proposed algorithm lends to a lower system cost  when the storage is reasonably efficient. From the figure, this corresponds to the case when the round-trip efficiency is greater than $0.65$, which includes the range of the round-trip efficiency for most energy storage in practice\cite{cg14}.   On the other hand, when the storage is highly inefficient, the greedy algorithm is shown to produce a better performance.

In Fig. \ref{fig:difumax}, we vary the value of the maximum charging and discharging rate $u_{i,\max}$ of all storage from $0.1$ kW to $3$ kW. Note that  for the greedy algorithm, $u_{i,\max}$ only  affects the constraints of the net charging and discharging amount, and for the proposed algorithm, $u_{i,\max}$ additionally affects the design of $V_{i,\max}$. We see that,  the system cost of the greedy algorithm decreases with $u_{i,\max}$, while the system cost of the proposed algorithm first decreases and then increases.  For the proposed algorithm, the increasing trend of the system cost  could be explained using Theorem \ref{the:per}.1, in which the  gap to the optimal objective value increases with $u_{i,\max}$. Moreover, from the figure, when $u_{i,\max}$ is less than $1.5$ kW, or, when the charging duration of the storage is larger than $6.6$  time units, the proposed algorithm outperforms the greedy one. Since the time scale we consider is seconds to minutes, this is the case for most batteries as the time scale of their charging duration is  hours\cite{cg14}. To improve the algorithm for large $u_{i,\max}$ is left for the  future.

\subsection{Effect of Other System Parameters}
In Fig. \ref{fig:power_flow}, we show the power flow $f_{i,t}$ between the substation and the $i$-th phase as well as their average, for $i=1, 2, 3$. Recall that the purpose of phase balancing is to make $f_{i,t}$ of all phases as close as possible. The figure shows that the curves of the power flows coincide most of the time. To further narrow the gap of these curves, we can increase the coefficient of the loss function $F(x)$ so as to impose more penalty for the flow deviation. In return, the system cost would be higher.

Although the three-phase transmission is dominant in practice, we are interested in finding  how the number of phases affects the algorithm performance. 
In Fig. \ref{fig:difn}, we increase the number of phases $N$ from $2$ to $8$. For both algorithms, the system cost grows linearly with $N$, which is expected since the  system cost sums up the costs of all phases. Moreover, as $N$ increases, the performance gain of the proposed algorithm over the greedy algorithm increases.

\subsection{Convergence of Distributed Implementation}
\begin{figure}[t]
	\centering
	\includegraphics[height= 2in,width = 2.5in]{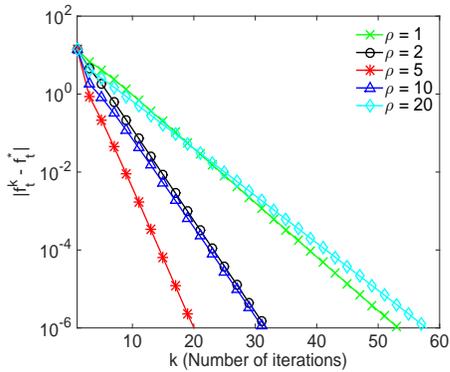}
	\caption{Performance gap versus number of iterations  for distributed algorithm.}
	\label{fig:admmrho}
\end{figure}
 In Fig. \ref{fig:admmrho}, we examine the convergence behavior of the distributed algorithm presented in Section III-B for various values of the $\rho$ parameter.  We show the gap between the  objective value of a per-slot optimization problem at iteration $k$ and its minimum objective value over iterations. We see that, for all $\rho$ values, the gap diminishes at a linear convergence rate. In particular, setting $\rho = 5$ leads to the best convergence performance. For a moderate accuracy requirement, the iterative procedure can be stopped within $20$ iterations. The fast convergence of the proposed algorithm is observed in general with appropriate $\rho$ values, and we omit the curves of other per-slot optimization problems for brevity.
 
 \vspace{1cm}
\section{Conclusion and Future Work}\label{sec:con}
We have investigated the problem of phase balancing with energy storage. 
We have proposed both centralized and distributed real-time algorithms for ideal energy storage and further extended the algorithms to accommodate non-ideal  energy storage. Moreover, we have conducted extensive simulation to evaluate the algorithm performance, showing that it can substantially outperform a greedy alternative. Our key conclusions are that positive correlations between the phases make phase balancing easier, and that evenly allocating storage over the phases results in the best performance.

For future work, we are interested in incorporating  system statistics into the algorithm design to further improve  performance, and also combining energy storage with traditional methods such as feeder reconfiguration for phase balancing.

\appendices
\section{Proof of Relaxation from \textbf{P1} to \textbf{P2}}\label{app:relax}
Using the  energy state update $s_{i,t+1} = s_{i,t} +u_{i,t}$,  we can derive  that the left hand side of  constraint \eqref{uitav} equals  the following:
\begin{align}\label{xst}
	\lim_{T\to\infty}\;\frac{1}{T}\sum_{t=0}^{T-1} \E[u_{i,t}] = \lim_{T\to\infty}\;\frac{\E[s_{i,T}]}{T}-\lim_{T\to\infty}\;\frac{\E[s_{i,0}]}{T}.
\end{align}
In \eqref{xst}, if $s_{i,t}$ is always bounded, \ie,  constraint \eqref{sit} holds, then the right hand side of  \eqref{xst} equals zero and thus constraint \eqref{uitav} is satisfied. Therefore, \textbf{P2} is  a relaxed problem of \textbf{P1}.

\section{An Upper Bound of the Drift-Plus-Cost Function}\label{app:lemdpp}
In the following lemma, we show that  the drift-plus-cost function is upper bounded.
\begin{lemma}\label{lem:dpp}
	For all possible decisions and all possible values of $s_{i,t}, i\in\mathcal{E}$, at each time slot $t$, the drift-plus-cost function is upper bounded as follows:
	\begin{align*}
		&\Delta(\bs_{t}) + \E[w_{t}|\bs_t]\\		
		\le & \E[w_{t}|\bs_t]+\sum_{i\in\mathcal{E}}\frac{u_{i,\max}^2}{2V_i} 
		+\frac{s_{i,t}-\beta_i}{V_i}\E\left[u_{i,t}|\bs_t\right].	
	\end{align*}
\end{lemma}
\begin{IEEEproof}
	Based on the definition of $L(s_{i,t})$ and the update of $s_{i,t}$,
	\begin{align*}
		& L(s_{i,t+1}) - L(s_{i,t}) \\
		=& \frac{1}{2}\left[(s_{i,t+1}-\beta_i)^2-(s_{i,t}-\beta_i)^2 \right]\\
		\le &(s_{i,t}-\beta_i)u_{i,t} + \frac{1}{2} u_{i,\max}^2.
	\end{align*}
	Using the  upper bound above for all phase $i\in\mathcal{E}$,   taking the conditional expectation, and then adding the term $\E[w_{t}|\bs_t]$  gives the desired upper bound.
\end{IEEEproof}

\section{Proof of Proposition \ref{pro:beta}}\label{app:probeta}
Since the per-slot problem \textbf{P3} includes all constraints of \textbf{P1} except the energy state constraint, the key of the feasibility proof is to show that the energy state $s_{i,t}$  is bounded within the interval $[s_{i,\min}, s_{i,\max}]$. To this end, we first prove the following lemma which gives a sufficient condition for charging or discharging.
\begin{lemma}\label{lem:sitcd}
	Under Algorithm \ref{alg:rtid}, for $i\in\mathcal{E}$,
	\begin{enumerate}
		\item if $s_{i,t}<\beta_i - V_i(p_{\max}+D'_{i,\max}+C'_{i,\max})$, then $u_{i,t}^* = u_{i,\max}$;
		\item if $s_{i,t}>\beta_i - V_i(p_{\min}+D'_{i,\min}+C'_{i,\min})$, then $u_{i,t}^* = -u_{i,\max}$.
	\end{enumerate}
\end{lemma}
\begin{IEEEproof}
	For simplicity of notation, we drop the time index $t$ in \textbf{P3}. Using constraint \eqref{bit} we replace $l_j$ with $u_j-f_j-r_j$ in the objective of \textbf{P3}. Next we solve \textbf{P3} through the partitioning method by first fixing the optimization variables $\bff$ and $u_j, j\neq i$, and then minimizing over $u_i$. The  optimization problem with respect to $u_i$ is as follows.
	\begin{align}
		\nonumber
		\min_{u_i}\;  pu_i + D_i(u_i)+C_i(u_i-f_i-r_i)+\frac{(s_i-\beta_i)u_i}{V_i},
		\st \eqref{uiti}.
	\end{align} 
	The derivative of the  objective above with respect to $u_i$ is $\frac{\partial(\cdot)}{\partial u_i} = p+D_i'(u_i)+C_i'(u_i-f_i-r_i)+\frac{(s_i-\beta_i)}{V_i}$. Therefore, if $s_i$ is upper bounded as shown in Lemma \ref{lem:sitcd}.1), we have $\frac{\partial(\cdot)}{\partial u_i}<0$ and thus $u_{i,t}^* =u_{i,\max}$. Or, if $s_i$ is lower bounded as shown  in Lemma \ref{lem:sitcd}.2), we have $\frac{\partial(\cdot)}{\partial u_i}>0$ and thus $u_{i,t}^* =-u_{i,\max}$. 
\end{IEEEproof}

Using Lemma \ref{lem:sitcd} above and the definition of $\beta_i$,  we can easily show the boundedness of the energy state  by mathematical induction, which is omitted here. 

\section{Proof of Theorem \ref{the:per}}\label{app:theper}
We prove Theorem \ref{the:per}.1 and Theorem \ref{the:per}.2 together. Denote $\tilde{w}$ as the optimal objective value of \textbf{P2}. In the following lemma, we show the existence of a special algorithm for \textbf{P2}. 
\begin{lemma}\label{lem:optst}
	For \textbf{P2}, there exists a stationary and randomized solution $\ba_t^s$ that only depends on the system state $\bq_t$, and at the same time satisfies the following conditions:
	\begin{align*}
		\E[w_t^s] \le \tilde{w}, \; \forall t,\quad
		\E[u_{i,t}^s] = 0, \; \forall i\in\mathcal{E}, t,
	\end{align*}
\end{lemma}
where  the expectations are taken over the randomness of the system state and the possible randomness of the actions.

The proof of Lemma \ref{lem:optst} follows  from Theorem 4.5 in \cite{bkneely} and  is omitted for brevity. 
Using Lemmas \ref{lem:dpp} and \ref{lem:optst}, the drift-plus-cost function under Algorithm \ref{alg:rtid} can be upper bounded as follows:
\begin{align}
	\nonumber
	&\Delta(\bs_{t}) + \E[w_{t}^*|\bs_t]\\
	\label{lyth0}
	\le & \E[w_t^s|\mb{s}_t] + \sum_{i\in\mathcal{E}}\Big[\frac{u_{i,\max}^2}{2V_i} + \frac{s_{i,t}-\beta_i}{V_i}\E\big[{u}_{i,t}^s|\mb{s}_t\big]\Big]\\
	\label{lyth1}
	\le &\tilde{w} + \sum_{i\in\mathcal{E}}\frac{u_{i,\max}^2}{2V_i}\\
	\label{lyth2}
	\le & w^{\opt} + \sum_{i\in\mathcal{E}}\frac{u_{i,\max}^2}{2V_i}
\end{align}
where  \eqref{lyth0} is derived based on Lemma \ref{lem:dpp} and the fact that \textbf{P3} minimizes the upper bound of the drift-plus-cost function,    \eqref{lyth1} is derived based on Lemma \ref{lem:optst} and the fact that the action $\ba_t^s$ is independent of $\mb{s}_t$, and the inequality in \eqref{lyth2} holds since \textbf{P2} is a relaxed problem of  \textbf{P1}. 

Taking  expectations over $\mb{s}_t$ on both sides of \eqref{lyth2} and summing over $t\in \{0,\cdots, T-1\}$ yields

\vspace{-0.4cm}
{\small{
		\begin{align}
			\nonumber
			\sum_{i\in\mathcal{E}}\E\Big[\frac{L(s_{i,T})-L(s_{i,0})}{V_i}\Big]+\sum_{t=0}^{T-1}\E[{w}_t^*]\le (w^{\opt} + \sum_{i\in\mathcal{E}}\frac{u_{i,\max}^2}{2V_i})T.
		\end{align}}}
		Note that $L(s_{i,T})$ is non-negative. Divide both sides of the above inequality by $T$.  After some arrangement,  there is
		\begin{align}\label{lydf2}
			\frac{1}{T}\sum_{t=0}^{T-1}\E[{w}_t^*] - w^{\opt} \le  \sum_{i\in\mathcal{E}}\Big[\frac{u_{i,\max}^2}{2V_i} + \frac{\E[L(s_{i,0})]}{TV_i} \Big],
		\end{align}
		which is the conclusion in Theorem \ref{the:per}.2.
		Taking $\limsup$ on both sides of \eqref{lydf2} gives Theorem \ref{the:per}.1.

		\section{Proof of Proposition \ref{pro:simax}}\label{app:prosimax}
		Denote $S$ as the fixed total energy capacity of storage. For simplicity of notation, we drop the index $i$ when the parameters are the same over all phases or storage units.
		Given the assumptions in Proposition \ref{pro:simax}, the optimization problem can be formulated as follows.
		\begin{align*}
			&\min_{s_{i,\max}}\; \textstyle{\sum_{i\in\mathcal{E}}\frac{u_{\max}^2(p_{\max}-p_{\min}+D'_{\max}-D'_{\min}+C'_{\max}-C'_{\min})}{2(s_{i,\max}-s_{\min}-2u_{\max})}}\\
			&\st \quad \sum_{i\in\mathcal{E}}s_{i,\max} = S
		\end{align*}
		where we have replaced $V_{i,\max}$ with its definition in \eqref{vimax}. It can be easily checked that the above problem is a convex optimization problem. Using the Karush-Kuhn-Tucker (KKT) conditions \cite{bkboyd}, the optimal solutions of $s_{i,\max}$ must be equal over $i$.
		
		\section{Proof of Theorem \ref{the:nonper}}\label{app:thenonper}
		1) To show the feasibility of $\{\ba^*\}$, it suffices to show that the resultant energy state $\hat{s}_{i,t}, i\in\mathcal{E},$ is bounded. First we give sufficient conditions of charging and discharging, which can be shown similarly to Lemma \ref{lem:sitcd}.
		\begin{lemma}\label{lem:sitcdnp}
			For $i\in\mathcal{E}$,
			\begin{enumerate}
				\item if $\hat{s}_{i,t}<\beta_i - V_i(\frac{p_{\max}}{\eta_i^+}+D'_{i,\max}+\frac{1}{\eta_i^+}C'_{i,\max})$, then $\hat{u}_{i,t}^+ = u_{i,\max}$;
				\item if $\hat{s}_{i,t}>\beta_i - V_i(p_{\min}\eta_i^-+D'_{i,\min}+\eta_i^-C'_{i,\min})$, then $\hat{u}_{i,t}^- = u_{i,\max}$.
			\end{enumerate}
		\end{lemma}
		Using Lemma \ref{lem:sitcdnp} and the mathematical induction arguments,
		we can show that $\hat{s}_{i,t}\in[s_{i,\min},s_{i,\max}], \forall i\in\mathcal{E}$. Note that the adjustment from  $(\hat{u}_i^{+},\hat{u}_i^{-})$ to $(\hat{u}_i^{+*},\hat{u}_i^{-*})$ does not change the difference  $\hat{u}_i^{+}-\hat{u}_i^{-}$. Therefore, the resultant energy state $s_{i,t}^*$ equals $\hat{s}_{i,t}$ and thus is bounded within $[s_{i,\min},s_{i,\max}]$.
		
		2) We prove Theorem \ref{the:nonper}.2 and Theorem \ref{the:nonper}.3 together.  Similar to the ideal case, the relaxed problem of \textbf{P1} can be formed as follows.
		\begin{align}
			\nonumber
			\textbf{P2':} \quad \min_{\{\ba_t\}}\quad &\limsup_{T\to\infty} \frac{1}{T}\sum_{t=0}^{T-1}\E[w_{t}]\\
			\nonumber
			\st \quad &\eqref{uit}, \eqref{uitn}, \eqref{fit}, \eqref{bit},\\
			\nonumber
			& \lim_{T\to\infty} \frac{1}{T}\sum_{t=0}^{T-1} \E[u_{i,t}^+-u_{i,t}^-] = 0, \forall i\in\mathcal{E}.
		\end{align}
		Denote the optimal value of \textbf{P2'} by $\tilde{w}'$. We first give the following two lemmas, which can be shown similarly to Lemmas \ref{lem:dpp} and \ref{lem:optst}.
		\begin{lemma}\label{lem:dppnp}
			For all possible decisions and all possible values of $s_{i,t}, i\in\mathcal{E}$, in each time slot $t$, the drift-plus-cost function is upper bounded as follows
			\begin{align}
				\nonumber
				&\Delta(\bs_{t}) + \E[w_{t}|\bs_t]\le \sum_{i\in\mathcal{E}}\frac{u_{i,\max}^2}{2V_i} +\E\big[w_{t}|\bs_t\big]\\
				\label{dppnpup}
				&\hspace{1cm}+\sum_{i\in\mathcal{E}}\frac{s_{i,t}-\beta_i}{V_i}\E\left[u^+_{i,t}-u^-_{i,t}|\bs_t\right].
			\end{align}
		\end{lemma}
		\begin{lemma}\label{lem:optstnp}
			For \textbf{P2'}, there exists a stationary and randomized solution $\ba_t^s$ that only depends on the system state $\bq_t$, and at the same time satisfies the following conditions:
			\begin{align}
				\label{npswt}
				&\E[w_t^s] \le \tilde{w}', \quad \forall t,\\
				\label{npuxt}
				&\E[u_{i,t}^{+s}-u_{i,t}^{-s}] = 0, \quad \forall i\in\mathcal{E}, t.
			\end{align}
		\end{lemma}
		
		Denote the optimal values of \textbf{P3'} under $\hat{\ba}_t$ and the adjusted solution $\ba^*_t$ by $\hat{g}_t$ and $g_t^*$, respectively. In the following lemma, we characterize the gap between $\hat{g}_t$ and $g_t^*$.
		\begin{lemma}\label{lem:gt}
			Under the proposed algorithm, at each time $t$ we have $g_t^* - \hat{g}_t\le \epsilon$, where $\epsilon\define \sum_{i\in\mathcal{E}}p_{\max}u_{i,\max}(\frac{1}{\eta_i^+}+\eta_i^-)+2D_{i,\max}+C_{i,\max}$.
		\end{lemma}
		\begin{IEEEproof}
			Using the objective of \textbf{P3'}, we have
			\begin{align}
				\nonumber
				&g_t^* - \hat{g}_t\\
				\nonumber
				&\le \sum_{i\in\mathcal{E}}p_t(\frac{1}{\eta_i^+}u_{i,t}^{+*}-\eta_i^-u_{i,t}^{-*}) + D_{i}(u_{i,t}^{+*}) +D_i(-u_{i,t}^{-*})+C_i(l^*_{i,t})\\
				\nonumber
				&-p_t(\frac{1}{\eta_i^+}\hat{u}_{i,t}^{+}-\eta_i^-\hat{u}_{i,t}^{-}) - D_{i}(\hat{u}_{i,t}^{+}) -D_i(-\hat{u}_{i,t}^{-})- C_i(\hat{l}_{i,t})\\
				\nonumber
				&\le\sum_{i\in\mathcal{E}}p_t\frac{1}{\eta_i^+}u_{i,t}^{+*} +p_t\eta_i^-\hat{u}_{i,t}^{-} + D_{i}(u_{i,t}^{+*}) +D_i(-u_{i,t}^{-*})
				+C_i(l^*_{i,t})\\
				\nonumber
				&\le \epsilon.
			\end{align}
		\end{IEEEproof}
		
		Using Lemmas \ref{lem:dppnp}, \ref{lem:optstnp}, and \ref{lem:gt}, the drift-plus-penalty function can be further upper bounded as follows.
		\begin{align}
			\nonumber
			&\Delta(\bs^*_{t}) + \E[w_{t}^*|\bs^*_t]\\
			\nonumber
			&\le\E\big[\hat{w}_{t}|\bs^*_t\big]+\sum_{i\in\mathcal{E}}\big[\frac{u_{i,\max}^2}{2V_i}+\frac{s^*_{i,t}-\beta_i}{V_i}\E\left[\hat{u}^+_{i,t}-\hat{u}^-_{i,t}|\bs^*_t\right]\big]+\epsilon\\
			\nonumber
			&\le  \E\big[{w}^s_{t}|\bs^*_t\big]+\sum_{i\in\mathcal{E}}\big[\frac{u_{i,\max}^2}{2V_i} +\frac{s^*_{i,t}-\beta_i}{V_i}\E\left[\hat{u}^{+s}_{i,t}-\hat{u}^{-s}_{i,t}|\bs^*_t\right]\big]+\epsilon\\
			\nonumber
			&\le\epsilon+\sum_{i\in\mathcal{E}}\frac{u_{i,\max}^2}{2V_i}+\tilde{w}'\\
			\nonumber
			&\le\epsilon+\sum_{i\in\mathcal{E}}\frac{u_{i,\max}^2}{2V_i}+w^{\opt}.
		\end{align}
		The remaining proof is similar to that for Theorem \ref{the:per} and is omitted for brevity.


\begin{IEEEbiography}[{\includegraphics[width=1in,height=1.25in,clip,keepaspectratio]{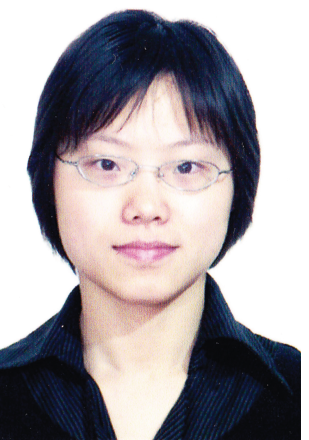}}]
	{Sun Sun} (S'11)
	received the B.S. degree in Electrical Engineering and Automation  from  Tongji University, Shanghai, China, in 2005. From 2006 to 2008, she was a software engineer in the Department of GSM Base Transceiver Station of Huawei Technologies Co. Ltd.. She received the M.Sc. degree in Electrical and Computer Engineering from University of Alberta, Edmonton, Canada, in 2011.
	Now, she is pursuing her Ph.D. degree in the Department of Electrical and Computer Engineering of University of Toronto, Toronto, Canada.
	Her current research interest lies in the areas of stochastic optimization, distributed control, learning, and economics, with the application of renewable generation, energy storage,
	demand response, and  power system operations.
\end{IEEEbiography}

\begin{IEEEbiography}[{\includegraphics[width=1in,height=1.25in,clip,keepaspectratio]{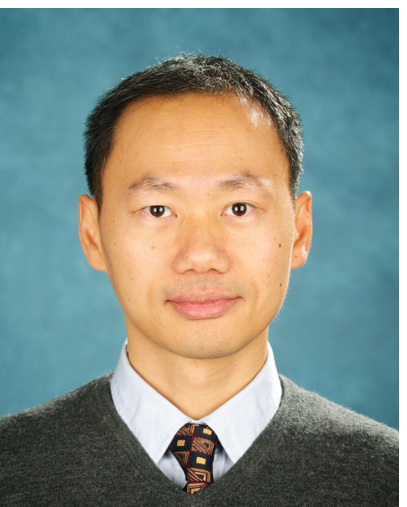}}]
	{Ben Liang}(S'94-M'01-SM'06) received honors-simultaneous B.Sc. (valedictorian) and M.Sc. degrees in Electrical Engineering from Polytechnic University in Brooklyn, New York, in 1997 and the Ph.D. degree in Electrical Engineering with a minor in Computer Science from Cornell University in Ithaca, New York, in 2001. In the 2001-2002 academic year, he was a visiting lecturer and post-doctoral research associate with Cornell University. He joined the Department of Electrical and Computer Engineering at the University of Toronto in 2002, where he is now a Professor. His current research interests are in mobile communications and networked systems. He has served as an editor for the IEEE Transactions on Communications, an editor for the IEEE Transactions on Wireless Communications, and an associate editor for the Wiley Security and Communication Networks journal, in addition to regularly serving on the organizational and technical committees of a number of conferences. He is a senior member of IEEE and a member of ACM and Tau Beta Pi.
\end{IEEEbiography}

\begin{IEEEbiography}[{\includegraphics[width=1in,height=1.25in,clip,keepaspectratio]{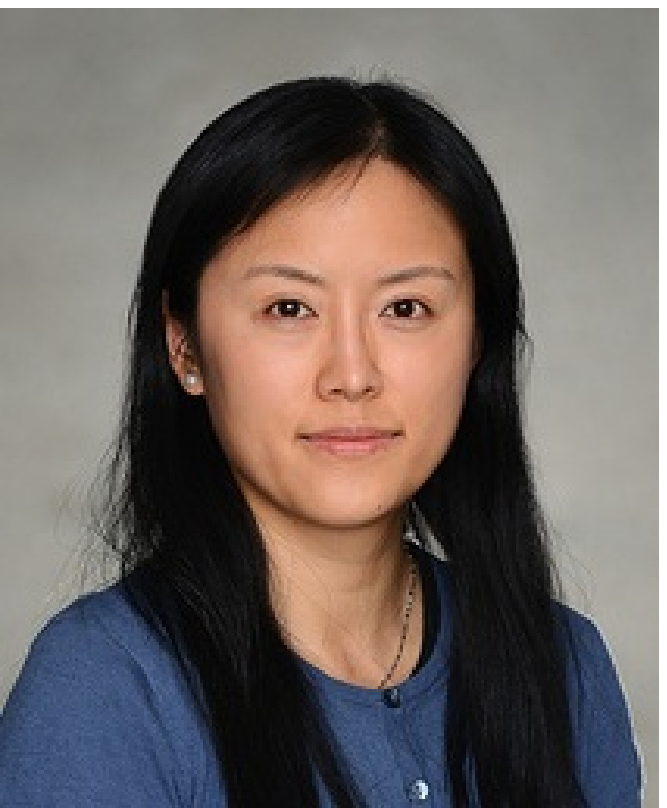}}]
	{Min Dong} (S'00-M'05-SM'09) received the B.Eng. degree from Tsinghua University, Beijing, China, in 1998, and the Ph.D. degree in electrical and computer engineering with minor in applied mathematics from Cornell University, Ithaca, NY, in 2004. From 2004 to 2008, she was with Corporate Research and Development, Qualcomm Inc., San Diego, CA. In 2008, she joined the Department of Electrical, Computer and Software Engineering at University of Ontario Institute of Technology, Ontario, Canada, where she is currently an Associate Professor. She also holds a status-only Associate Professor appointment with the Department of Electrical and Computer Engineering, University of Toronto since 2009. Her research interests are in the areas of statistical signal processing for communication networks, cooperative communications and networking techniques, and stochastic network optimization in dynamic networks and systems.
	
	Dr. Dong received the Early Researcher Award from Ontario Ministry of Research and Innovation in 2012, the Best Paper Award at IEEE ICCC in 2012, and the 2004 IEEE Signal Processing Society Best Paper Award. She served as an Associate Editor for the IEEE TRANSACTIONS ON SIGNAL PROCESSING from 2010 to 2014, and as  an Associate Editor for the IEEE SIGNAL PROCESSING LETTERS from 2009 to 2013. She was a Symposium lead co-chair of the Communications and Networks to Enable the Smart Grid Symposium at the IEEE International Conference on Smart Grid Communications (SmartGridComm) in 2014.
	She has been an elected member of IEEE Signal Processing Society Signal Processing for Communications and Networking (SP-COM) Technical Committee since 2013. 
\end{IEEEbiography}

\begin{IEEEbiography}[{\includegraphics[width=1in,height=1.25in,clip,keepaspectratio]{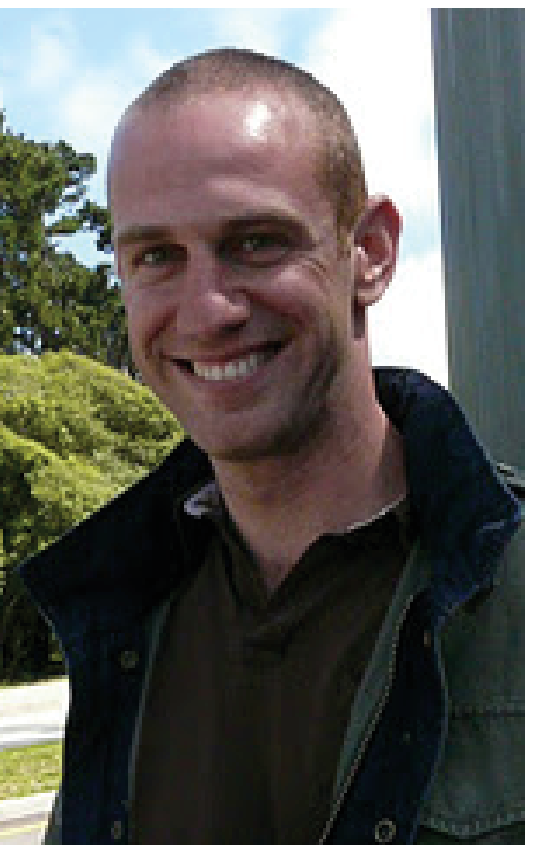}}]
	{Joshua A. Taylor} (M'11) received the B.S. degree from Carnegie Mellon University in 2006, and the S.M. and Ph.D. degrees from the Massachusetts Institute of Technology in 2008 and 2011, all in Mechanical Engineering. From 2011 to 2012, he was a postdoctoral researcher in Electrical Engineering and Computer Sciences at the University of California, Berkeley. He is currently an assistant professor in the Department of Electrical and Computer Engineering and the associate director of the Institute for Sustainable Energy at the University of Toronto. His current research interests include renewable energy and demand response, machine learning, and infrastructural couplings.
\end{IEEEbiography}

\end{document}